\documentclass[11pt,a4paper]{article}

\usepackage[T1]{fontenc}
\usepackage[utf8]{inputenc}
\usepackage{amsmath,amssymb,amsthm}
\usepackage{enumitem}
\usepackage[left=0.2cm, right=0.2cm, top=1.1cm, bottom=1.1cm]{geometry}

\usepackage{xcolor}
\definecolor{linkblue}{rgb}{0.1,0.2,0.55}
\usepackage[colorlinks=true,linkcolor=linkblue,citecolor=linkblue,urlcolor=linkblue]{hyperref}

\newtheorem{theorem}{Theorem}[section]
\newtheorem{proposition}[theorem]{Proposition}
\newtheorem{lemma}[theorem]{Lemma}
\newtheorem{corollary}[theorem]{Corollary}
\newtheorem{conjecture}[theorem]{Conjecture}
\theoremstyle{definition}
\newtheorem{definition}[theorem]{Definition}
\newtheorem{example}[theorem]{Example}
\newtheorem{remark}[theorem]{Remark}

\numberwithin{equation}{section}
\DeclareMathOperator{\des}{des}
\DeclareMathOperator{\inv}{inv}

\title{Beyond geometric lattices: pseudo-atoms and redundant join representations\\ in the face lattice of the permutohedron}
\author{\Large{Giovanni Rossi}\\ Independent researcher, Italy\\ \small{\textit{Email:} \texttt{roxyjean@gmail.com} - ORCID \textit{identifier:} 0000-0001-5215-4325}}
\date{\today}

\begin{document}
\maketitle
\begin{abstract}
We study the lattice $\mathcal L_n$ obtained from the poset $\mathrm{OP}_n$ of ordered set partitions of $[n]$ (the Fubini poset), ordered by merging adjacent blocks, by adjoining a
bottom element $\bot$ of rank $-1$. It is the face lattice of the permutohedron. We look at it from the point of view of geometric lattices: it is graded, atomic and coatomic in the strong sense that every element is a meet of coatoms, but it is not semimodular. Its atoms are the $n!$ permutations, and its elements of rank one, which we call pseudo-atoms, are the edges of the permutohedron. Describing ordered set partitions by chains of subsets, the cuts, we show that the rank is modular whenever the meet is above $\bot$, while the join of two distant atoms can have large rank, so that the rank of the join is a metric on permutations. An element has a unique representation as a join of atoms only up to rank one; above it, irredundant representations have sizes ranging from two to the rank, and the redundancy of an element is its number of atoms minus two. The M\"obius function is that of an Eulerian lattice, the Whitney numbers are the face numbers of the permutohedron, and the characteristic polynomial has the alternating coefficients of the geometric case, but only $\bot$ and $\top$ are left-modular, so the lattice is not supersolvable. We show that it is CL-shellable, we exhibit EL-labelings for $n\le5$, and we compare the independent sets of Raz \cite[arXiv:2307.05005 (v3, Jan. 2026)]{Raz} with irredundant families and with the matroid of the vertices of the permutohedron. Replacing ordered set partitions by partial flags of subspaces of $\mathbb F_q^n$ gives a $q$-analogue, whose M\"obius function is $(-1)^{\mathrm{rk}+1}q^{\kappa}$, with $\kappa$ the number of pairs of elements in a common block at $q=1$, and whose characteristic polynomial and Betti numbers are polynomials in $q$. Open problems are listed.
\end{abstract}

\noindent\textbf{Keywords.} ordered set partitions, Fubini poset, permutohedron, face lattice, geometric lattice, semimodularity, join representations, M\"obius function, characteristic polynomial, supersolvable lattice, CL-shellability.

\noindent\textbf{MSC 2020.} Primary 06A07; secondary 05A18, 05B35, 05E45, 06C10, 52B05.

\section{Introduction}\label{sec:intro}

Geometric lattices---the finite atomic semimodular lattices, equivalently the lattices of flats of matroids---form one of the best understood classes of graded posets \cite{StanleyEC1,Oxley}. Their rank function, M\"obius function and Whitney numbers are governed by matroid theory, and semimodularity is what makes that theory work. Many lattices that arise in combinatorics are nevertheless graded, atomic and even coatomic without being semimodular; face lattices of polytopes are the standard examples \cite{Ziegler}. In this paper we study one such lattice, which is simple to define, for which many structural questions can be answered in closed form, and for which some remain open.

\medskip\noindent\textbf{The object.}
Let $\mathrm{OP}_n$ be the set of ordered set partitions
$x=(B_1\mid\cdots\mid B_k)$ of $[n]=\{1,\dots,n\}$, ordered by $x\le y$ when $y$
is obtained from $x$ by merging adjacent blocks. Its cardinality is the Fubini
number $\sum_k k!\,S(n,k)$ \cite{StanleyEC1}; the minimal elements are the $n!$
permutations (all blocks singletons) and the maximum is $([n])$. Any two
elements of $\mathrm{OP}_n$ have a join, but not a meet. We adjoin a bottom
element $\bot$, declare $\mathrm{rk}(\bot)=-1$, and put
$\mathcal L_n=\mathrm{OP}_n\cup\{\bot\}$; this is a lattice
(Theorem~\ref{thm:lattice}), with maximum $\top=([n])$. The permutations, which
we call \emph{atoms}, have rank $0$.

The lattice $\mathcal L_n$ is the face lattice of the permutohedron
\cite{Ziegler,Postnikov}, \cite[Proposition~1.4]{BilleraSarangarajan}: the element $(B_1\mid\cdots\mid B_k)$ is the face
whose vertices are the permutations obtained by ordering each block, and $\bot$
is the empty face. Its dual is the face poset of the braid arrangement (the
Coxeter complex of type $A_{n-1}$) with a bottom element adjoined \cite{BjornerBrenti,AbramenkoBrown}.
It was studied by Sagan and Sundaram \cite{SaganSundaram}, who prove that it is
atomic, not semimodular in general, with Boolean upper intervals and M\"obius
function $(-1)^{\mathrm{rk}\,y-\mathrm{rk}\,x}$ (Lemma~2.1 and Theorem~2.4 there), and who give a recursive atom
ordering (Theorem~2.7), which yields CL-shellability; they ask whether the family of lattices containing $\mathcal L_n$ (the case $d=1$) has EL-labelings (Question~2.8). Face lattices of polytopes are CL-shellable in general
\cite[Theorem~4.5]{BjornerWachs}. We make no claim of novelty for these facts. Our interest
lies elsewhere: in how much of the geometric-lattice toolkit survives for
$\mathcal L_n$, and in what replaces it when it does not.

\medskip\noindent\textbf{Pseudo-geometric lattices.}
We write $L$ for a generic finite lattice. We call $L$ \emph{pseudo-geometric} if
it is graded with $\mathrm{rk}(\bot)=-1$, atomic, and coatomic in the strong
sense that every element is a meet of coatoms. Geometric lattices (after shifting
the rank by one) and face lattices of polytopes are pseudo-geometric, so the class
is broad and we do not attempt to classify it. What it isolates is the loss of
semimodularity. In $\mathcal L_n$ the elements of rank $1$, which we call
\emph{pseudo-atoms} (they are the edges of the permutohedron), are the first
level above the atoms, and joining an element with an atom may raise its rank by
much more than one.

\medskip\noindent\textbf{Redundant join representations.}
In a geometric lattice, every family of atoms that is minimal with a given join $x$
has the same size, the rank of $x$ counted from the atoms. In $\mathcal L_n$ this
fails in a precise way. An element $x$ lies above $s(x)=\prod_j|B_j|!$ atoms, and its
representation as a join of atoms is unique if and only if $\mathrm{rk}\,x\le1$; the
pseudo-atoms are the elements of rank $1$, each the join of exactly two atoms. For
$\mathrm{rk}\,x\ge2$ two atoms already suffice to generate $x$, minimal families have
sizes ranging from $2$ to $\mathrm{rk}\,x$, and the closure operator
$S\mapsto\mathrm{At}(\bigvee S)$ does not satisfy the exchange property. The
resulting redundancy $s(x)-2$ is measured on the atoms by the rank of the join,
$\delta(\rho,\rho')=\mathrm{rk}(\rho\vee\rho')$, which is a metric on permutations.

\medskip\noindent\textbf{Cuts.}
The key device is the \emph{cut model}. To $x=(B_1\mid\cdots\mid B_k)$ we
associate the chain $C(x)=\{A_1\subset\cdots\subset A_{k-1}\}$, where $A_i$ is
the union of the first $i$ blocks of $x$. The map $x\mapsto C(x)$ is a bijection
from $\mathrm{OP}_n$ onto the chains of proper nonempty subsets of $[n]$, with
$x\le y$ if and only if $C(y)\subseteq C(x)$, and $\mathrm{rk}(x)=n-1-|C(x)|$.
Joins become intersections of chains, and meets become unions when the union is
again a chain and equal $\bot$ otherwise. The same description of the faces of the permutohedron as chains of subsets is in
\cite[Proposition~1.3 and the paragraph after it]{BilleraSarangarajan}. Everything below is a consequence of
this dictionary.

\medskip\noindent\textbf{Related work.}
Besides \cite{SaganSundaram}, the closest framework is that of Raz \cite{Raz},
who studies independent sets of finite atomic graded lattices that are not
geometric, defined through non-broken-circuit sets for orderings of the atoms in
the sense of Blass and Sagan \cite{BlassSagan}. In Section~\ref{sec:open} we compare
them with our irredundant families, and we show that the independent sets of
$\mathcal L_n$ do not form a matroid for $n\ge4$, although $\mathcal L_n$ embeds in
the lattice of flats of the matroid of the vertices of the permutohedron. The
redundancy of join representations has a counterpart in the theory of frames
\cite{BCHL}, to which we return in a remark.

\medskip\noindent\textbf{Plan of the paper.}
Throughout, $n\ge 3$ unless stated otherwise. The paper proceeds in the
following steps.
\begin{enumerate}
\item (Section~\ref{sec:cut}) We establish the cut model and show that
  $\mathcal L_n$ is graded, with $\mathrm{rk}(x)=n-1-|C(x)|$.
\item (Section~\ref{sec:lattice}) We express joins and meets in cut terms, show
  that every upper interval $[x,\top]$ is Boolean, and deduce that the lattice
  is atomic and coatomic in the strong sense, hence pseudo-geometric. We also show that
  the rank is modular,
  $\mathrm{rk}\,x+\mathrm{rk}\,y=\mathrm{rk}(x\wedge y)+\mathrm{rk}(x\vee y)$,
  whenever $x\wedge y>\bot$. When the meet is $\bot$ the defect can have
  either sign for $n\ge 3$, so the lattice is not semimodular; in particular the
  join of an element with an atom need not cover it
  (Subsection~\ref{sec:semimod}).
\item (Section~\ref{sec:pseudo}) We describe the atoms (the permutations), the
  $2^n-2$ coatoms and the pseudo-atoms, show that every element of rank at
  least $1$ is a join of pseudo-atoms, and that the representation of an element
  as a join of atoms is unique if and only if its rank is at most $1$.
\item (Section~\ref{sec:metric}) We study the size $s(x)$, the closure operator on
  atoms and its failure of exchange, and the sizes of irredundant families; we
  show that $\delta(\rho,\rho')=\mathrm{rk}(\rho\vee\rho')$ is a metric on atoms,
  bounded by Kendall's $\tau$ and by $n-1$, and interpret it through cuts and
  through the inversion graph: far-apart atoms have joins of large rank.
\item (Section~\ref{sec:mobius}) We prove $\mu(\bot,x)=(-1)^{\mathrm{rk}(x)+1}$,
  so the lattice is Eulerian, compute the Whitney numbers
  $(n-k)!\,S(n,n-k)$, give the characteristic polynomial, and conjecture that its
  roots are real.
\item (Section~\ref{sec:modular}) For $n\ge 3$ no element other than $\bot$ and
  $\top$ is left-modular, hence the lattice is not supersolvable \cite{Stanley72,BlassSagan}.
\item (Section~\ref{sec:shell}) We show that $\mathcal L_n$ is CL-shellable, with
  every linear extension of the weak order on permutations giving a recursive atom ordering,
  discuss the flag $h$-vector, and discuss EL-labelings, which we found by computer for
  $n\le5$; this bears on a question of Sagan and Sundaram
  \cite[Question~2.8]{SaganSundaram}.
\item (Section~\ref{sec:qanalogue}) We replace ordered set partitions by partial flags of subspaces of
  $\mathbb F_q^n$. The M\"obius function becomes $(-1)^{\mathrm{rk}+1}q^{\kappa}$, by the Solomon--Tits theorem,
  and the characteristic polynomial and the numbers $\beta_j$ of Theorem~\ref{thm:chi} become polynomials in $q$.
\item (Section~\ref{sec:open}) We compare $\mathcal L_n$ with geometric lattices,
  prove that the rank is subadditive exactly when the meet is above $\bot$,
  relate $\mathcal L_n$ to the vertex matroid of the permutohedron and to the
  independent sets of Raz, and collect open problems.
\end{enumerate}

\section{The cut model}\label{sec:cut}

Throughout, $n\ge 3$ unless stated otherwise and $[n]=\{1,\dots,n\}$. A subset
$A\subseteq[n]$ is \emph{proper} if $\emptyset\ne A\ne[n]$.

\subsection{Ordered set partitions}

An \emph{ordered set partition} of $[n]$ is a sequence
$x=(B_1\mid\cdots\mid B_k)$ of nonempty, pairwise disjoint sets, called
\emph{blocks}, whose union is $[n]$. We write $k(x)=k$ for the number of blocks
and $\mathrm{OP}_n$ for the set of all ordered set partitions of $[n]$. For
$x=(B_1\mid\cdots\mid B_k)$ and $y\in\mathrm{OP}_n$ we write $x\le y$ if $y$
is obtained from $x$ by merging adjacent blocks, that is, if there are indices
$0=i_0<i_1<\cdots<i_l=k$ such that $y=(D_1\mid\cdots\mid D_l)$ with
\[
D_j=\bigcup_{i_{j-1}<m\le i_j}B_m\qquad(1\le j\le l).
\]
Thus each block of $y$ is the union of a collection of consecutive blocks of $x$ (at
least one block, so that $\le$ is reflexive), and these collections, taken in
order, exhaust
$B_1,\dots,B_k$. Merges may occur anywhere in the sequence, as long as they are consecutive.
This relation is reflexive and transitive, and it is antisymmetric because a
proper merging strictly decreases the number of blocks; so $\le$ is a partial
order (it also follows from Proposition~\ref{prop:dictionary} below). The element
$\top=([n])$, with a single block, is the maximum of $\mathrm{OP}_n$. The
cardinality of $\mathrm{OP}_n$ is the Fubini number
$\sum_{k\ge 1}k!\,S(n,k)$, with $S(n,k)$ denoting the Stirling numbers of the second kind \cite{StanleyEC1}.

\begin{definition}[Cuts]\label{def:cuts}
Let $x=(B_1\mid\cdots\mid B_k)\in\mathrm{OP}_n$. For $1\le i\le k-1$ put
$A_i(x)=B_1\cup\cdots\cup B_i$, the union of the first $i$ blocks of $x$; we also put $A_0(x)=\emptyset$ and
$A_k(x)=[n]$. Let
\[
C(x)=\{A_1(x),\dots,A_{k-1}(x)\}.
\]
Then $A_1(x)\subset\cdots\subset A_{k-1}(x)$ is a chain of proper subsets of
$[n]$, so $|C(x)|=k(x)-1$. We call $C(x)$ the \emph{cut set} of $x$ and its
elements the \emph{cuts} of $x$. Let $\mathcal C_n$ be the set of all chains of
proper subsets of $[n]$, the empty chain included.
\end{definition}

\begin{proposition}[The cut dictionary]\label{prop:dictionary}
The map $x\mapsto C(x)$ is a bijection from $\mathrm{OP}_n$ onto
$\mathcal C_n$, and for $x,y\in\mathrm{OP}_n$
\[
x\le y\quad\Longleftrightarrow\quad C(y)\subseteq C(x).
\]
\end{proposition}

\begin{proof}
The inverse of $C$ sends a chain $A_1\subset\cdots\subset A_k$ to the ordered set
partition $(A_1\mid A_2\setminus A_1\mid\cdots\mid[n]\setminus A_k)$, whose
blocks are nonempty because the inclusions are strict and the sets are proper.
So $C$ is a bijection.

If $y$ is obtained from $x$ by merging adjacent blocks, every cut of
$y$ is a cut of $x$ (indeed $A_j(y)=A_{i_j}(x)$), hence
$C(y)\subseteq C(x)$. Conversely, let
$C(y)\subseteq C(x)$ and write $y=(D_1\mid\cdots\mid D_l)$, with $A_0(y)=\emptyset$
and $A_l(y)=[n]$. For each $j$ the sets $A_{j-1}(y)\subset A_j(y)$ are
cuts of $x$ or one of $\emptyset$, $[n]$ (that is, of the form $A_q(x)$ with
$0\le q\le k$), say
$A_{j-1}(y)=A_q(x)$ and $A_j(y)=A_p(x)$ with $q<p$. Then
$D_j=A_j(y)\setminus A_{j-1}(y)=B_{q+1}\cup\cdots\cup B_p$ is a union of
consecutive blocks of $x$, and these unions follow one another for
$j=1,\dots,l$. Hence $y$ is obtained from $x$ by merging adjacent blocks.
\end{proof}

Thus going up in $\mathrm{OP}_n$ means \emph{deleting cuts}: the fewer cuts, the
coarser the partition. The permutations, which have the maximum number $n-1$ of
cuts, are at the bottom, and $\top$, which has none, is at the top.

\begin{definition}[Rank and the poset $\mathcal L_n$]\label{def:rank}
For $x\in\mathrm{OP}_n$ put
\[
\mathrm{rk}(x)=n-k(x)=n-1-|C(x)|.
\]
Adjoin to $\mathrm{OP}_n$ a new element $\bot$ below every element, set
$\mathrm{rk}(\bot)=-1$, and put $\mathcal L_n=\mathrm{OP}_n\cup\{\bot\}$.
\end{definition}

We show in Section~\ref{sec:lattice} that $\mathcal L_n$ is a lattice
(Theorem~\ref{thm:lattice}); in this section we only need the order and the
rank. The elements of rank $0$ are the \emph{atoms} of $\mathcal L_n$, and those
of rank $1$ the \emph{pseudo-atoms}; the names are justified by
Corollary~\ref{cor:graded} and are studied in Section~\ref{sec:pseudo}.

\subsection{Covers, atoms and gradedness}

\begin{lemma}[Covers]\label{lem:covers}
Let $x,y\in\mathrm{OP}_n$ with $x<y$. The following are equivalent.
\begin{enumerate}
\item[(a)] $y$ covers $x$ in $\mathrm{OP}_n$.
\item[(b)] $C(y)=C(x)\setminus\{A\}$ for some $A\in C(x)$.
\item[(c)] $y$ is obtained from $x$ by merging two adjacent blocks.
\end{enumerate}
In this case $\mathrm{rk}(y)=\mathrm{rk}(x)+1$. Moreover every $x\ne\top$ is
covered by some element of $\mathrm{OP}_n$.
\end{lemma}

\begin{proof}
(b)$\Leftrightarrow$(c): merging the blocks $B_i$ and $B_{i+1}$ of $x$ deletes
exactly the cut $A_i(x)$ and keeps the others.

(a)$\Rightarrow$(b): since $x<y$ we have $C(y)\subsetneq C(x)$ by
Proposition~\ref{prop:dictionary}. Pick $A\in C(x)\setminus C(y)$, and let $z$ be
the element with $C(z)=C(x)\setminus\{A\}$. Then $x<z\le y$, so $z=y$ if $y$
covers $x$.

(b)$\Rightarrow$(a): if $x\le z\le y$ then
$C(x)\setminus\{A\}=C(y)\subseteq C(z)\subseteq C(x)$, so $C(z)$ equals one of
the two ends, that is, $z=y$ or $z=x$.

The rank statement follows from $|C(y)|=|C(x)|-1$. Finally, if $x\ne\top$ then
$C(x)\ne\emptyset$, and deleting any cut gives an element covering $x$.
\end{proof}

For a permutation $\rho=(\rho_1\mid\cdots\mid\rho_n)$, regarded as an ordered set
partition into singletons, the cut set is the complete flag
$F(\rho)=\{\{\rho_1\},\{\rho_1,\rho_2\},\dots,\{\rho_1,\dots,\rho_{n-1}\}\}$.
By Proposition~\ref{prop:dictionary}, $\rho\mapsto F(\rho)$ is a bijection from
the $n!$ permutations onto the maximal chains of proper subsets of $[n]$, and we
call $F(\rho)$ the \emph{flag} of $\rho$.

\begin{lemma}[Atoms below an element]\label{lem:atoms}
The minimal elements of $\mathrm{OP}_n$ are exactly the $n!$ permutations, and
they are the elements of rank $0$. For $x=(B_1\mid\cdots\mid B_k)$ a permutation
$\rho$ satisfies $\rho\le x$ if and only if $F(\rho)\supseteq C(x)$, and the
number of such permutations is $\prod_{j=1}^k|B_j|!$. They are the concatenations
of orderings of the blocks $B_1,\dots,B_k$.
\end{lemma}

\begin{proof}
An element with a block of size at least $2$ is not minimal, since splitting that
block in two parts gives a smaller element (Lemma~\ref{lem:covers}); a
permutation is minimal, because $x<y$ implies $k(y)<k(x)$, so the rank strictly
increases along the order and no element has rank below $0$. Hence the minimal
elements are the permutations. The
equivalence $\rho\le x\Leftrightarrow F(\rho)\supseteq C(x)$ is
Proposition~\ref{prop:dictionary}. A maximal chain $F\supseteq C(x)$ consists of
$C(x)$ together with, for each $j$, a maximal chain of subsets strictly between
$A_{j-1}(x)$ and $A_j(x)$ (with $A_0(x)=\emptyset$ and $A_k(x)=[n]$); the
latter amounts to an ordering of $B_j$, and there are $|B_j|!$ of them.
\end{proof}

\begin{corollary}[Gradedness]\label{cor:graded}
The poset $\mathcal L_n$ is graded, with rank function $\mathrm{rk}$: its
minimum $\bot$ has rank $-1$, it is covered exactly by the $n!$ permutations,
which have rank $0$, every cover raises the rank by one, and its maximum $\top$
has rank $n-1$. In particular every maximal chain of $\mathcal L_n$ has $n+1$
elements.
\end{corollary}

\begin{proof}
The covers within $\mathrm{OP}_n$ raise the rank by one by
Lemma~\ref{lem:covers}. The elements covering $\bot$ are the minimal elements of
$\mathrm{OP}_n$, which are the permutations by Lemma~\ref{lem:atoms}. By
Lemma~\ref{lem:covers} every element other than $\top$ is covered by another one,
so every maximal chain ends at $\top$, and its ranks run from $-1$ to $n-1$.
\end{proof}

\begin{example}\label{ex:n3}
Write subsets of $[3]$ as words, for instance $12=\{1,2\}$. The $13$ elements of
$\mathrm{OP}_3$, with their cut sets, are:
\begin{center}
\begin{tabular}{@{}cll@{}}
$\mathrm{rk}$ & elements $x$ & cut sets $C(x)$\\\hline
$0$ & $(1|2|3),\ (1|3|2),\ (2|1|3)$ & $\{1,12\},\ \{1,13\},\ \{2,12\}$\\
$0$ & $(2|3|1),\ (3|1|2),\ (3|2|1)$ & $\{2,23\},\ \{3,13\},\ \{3,23\}$\\
$1$ & $(12|3),\ (13|2),\ (23|1)$ & $\{12\},\ \{13\},\ \{23\}$\\
$1$ & $(1|23),\ (2|13),\ (3|12)$ & $\{1\},\ \{2\},\ \{3\}$\\
$2$ & $(123)=\top$ & $\emptyset$
\end{tabular}
\end{center}
So $(1|2|3)\le(12|3)$ and $(1|2|3)\le(1|23)$, since $\{12\}$ and $\{1\}$ are
cuts of $(1|2|3)$, while $(1|2|3)\not\le(13|2)$. For $n=4$,
$x=(24\mid1\mid3)$ has $C(x)=\{24,124\}$ and $\mathrm{rk}(x)=4-3=1=n-1-|C(x)|$.
\end{example}

\begin{remark}[Geometry]\label{rem:geometry}
By Lemma~\ref{lem:atoms}, the permutations below $x$ number
$\prod_j|B_j|!$. This is the number of vertices of the face of the
permutohedron which is a product of permutohedra on the blocks, a face of
dimension $\sum_j(|B_j|-1)=\mathrm{rk}(x)$; the face lattice of the permutohedron
is recovered in this way \cite[Proposition~1.4]{BilleraSarangarajan}. We do not use this identification in
Sections~\ref{sec:cut}--\ref{sec:metric}: all arguments there go through the cut model. It is used in
Sections~\ref{sec:shell} and~\ref{sec:open} (Proposition~\ref{prop:dictionary}).
\end{remark}

\section{Joins, meets and pseudo-geometricity}\label{sec:lattice}

We use the cut model of Section~\ref{sec:cut}: by Proposition~\ref{prop:dictionary}
an element $x\in\mathrm{OP}_n$ is the same thing as its cut set $C(x)$, and $x\le y$
means $C(y)\subseteq C(x)$. All the arguments below are statements about chains of
subsets of $[n]$.

\subsection{Joins and meets}

\begin{theorem}[Joins and meets]\label{thm:lattice}
Let $x,y\in\mathrm{OP}_n$.
\begin{enumerate}
\item The join $x\vee y$ exists in $\mathcal L_n$, it lies in $\mathrm{OP}_n$, and
  $C(x\vee y)=C(x)\cap C(y)$.
\item If $C(x)\cup C(y)$ is a chain, then $x\wedge y$ exists in $\mathrm{OP}_n$ and
  $C(x\wedge y)=C(x)\cup C(y)$. Otherwise the only common lower bound of $x$ and
  $y$ in $\mathcal L_n$ is $\bot$, so $x\wedge y=\bot$.
\item $\mathcal L_n$ is a lattice, with minimum $\bot$ and maximum $\top$. It is the
  smallest lattice containing $\mathrm{OP}_n$: if $L$ is a lattice containing
  $\mathrm{OP}_n$ as a subposet, the inclusion of $\mathrm{OP}_n$ extends to an order
  embedding of $\mathcal L_n$ into $L$.
\end{enumerate}
\end{theorem}

\begin{proof}
(1) The element $\bot$ is not an upper bound of $x$, so the upper bounds of $x$ and
$y$ in $\mathcal L_n$ are those in $\mathrm{OP}_n$, and there is at least one, namely
$\top$. By Proposition~\ref{prop:dictionary}, $z$ is an upper bound of $x$ and $y$ if
and only if $C(z)\subseteq C(x)\cap C(y)$. The set $C(x)\cap C(y)$ is contained in a
chain, hence is a chain, so it equals $C(z_0)$ for a unique $z_0\in\mathrm{OP}_n$, and
every upper bound $z$ satisfies $z\ge z_0$. Thus $z_0=x\vee y$.

(2) A lower bound $w\in\mathrm{OP}_n$ of $x$ and $y$ satisfies
$C(w)\supseteq C(x)\cup C(y)$. Such a $w$ exists only if $C(x)\cup C(y)$ is a chain,
and then the smallest possible $C(w)$ is the union itself, which gives the greatest
lower bound. If the union is not a chain, no element of $\mathrm{OP}_n$ lies below both
$x$ and $y$, and $\bot$ is the only common lower bound.

(3) By (1), any two elements of $\mathrm{OP}_n$ have a join; joins with $\bot$ are
immediate ($x\vee\bot=x$). A finite poset with a minimum and in which any two elements
have a join is a lattice, and the meets are those in (2) together with $x\wedge\bot=\bot$.
The maximum is $\top=([n])$ and the minimum is $\bot$.

For the last statement, let $L$ be a lattice containing $\mathrm{OP}_n$ as a subposet
and let $m$ be the meet in $L$ of the $n!$ permutations. Every element of
$\mathrm{OP}_n$ lies above some permutation, so $m\le x$ for all
$x\in\mathrm{OP}_n$. Moreover $m\notin\mathrm{OP}_n$: otherwise $m$ would be a common
lower bound in $\mathrm{OP}_n$ of $(1\mid 2,\dots,n)$ and $(2\mid 1,3,\dots,n)$, which
have none (Remark~\ref{rem:nomeets}). Send $\bot$ to $m$ and every $x\in\mathrm{OP}_n$ to
itself. The order is preserved, and also reflected: $m\le x$ for every $x$, while
$x\le m$ is impossible, since $x$ would then lie below every permutation, whereas a
permutation is minimal in $\mathrm{OP}_n$ and there are at least two of them.
\end{proof}

\begin{remark}\label{rem:nomeets}
Inside $\mathrm{OP}_n$, any two elements have a join, but they need not have a meet:
the coatoms $(1\mid 2,\dots,n)$ and $(2\mid 1,3,\dots,n)$ have cut sets $\{\{1\}\}$ and
$\{\{2\}\}$ whose union is not a chain, so they have no common lower bound at all. This
is the reason for adjoining $\bot$, and by Theorem~\ref{thm:lattice}(3) this is the
smallest way to obtain a lattice.
\end{remark}

\begin{example}\label{ex:joinmeet}
Let $n=4$ and write subsets as words. The atoms $x=(1|2|3|4)$ and $y=(1|3|2|4)$ have
$C(x)=\{1,12,123\}$ and $C(y)=\{1,13,123\}$, so $C(x)\cap C(y)=\{1,123\}$ and
$x\vee y=(1|23|4)$. The elements $(1|23|4)$ and $(1|2|34)$ have cut sets $\{1,123\}$
and $\{1,12\}$, whose union $\{1,12,123\}$ is a chain, so their meet is the atom
$(1|2|3|4)$. By contrast $(12|3|4)$ and $(13|2|4)$ have cut sets $\{12,123\}$ and
$\{13,123\}$, whose union is not a chain, so their meet is $\bot$.
\end{example}

\subsection{Upper intervals are Boolean}

\begin{proposition}[Upper intervals]\label{prop:boolean}
For $x\in\mathrm{OP}_n$ the map $y\mapsto C(x)\setminus C(y)$ is an isomorphism from
the interval $[x,\top]$ onto the Boolean lattice of subsets of $C(x)$, ordered by
inclusion. In particular $y$ covers $x$ if and only if $C(y)=C(x)\setminus\{A\}$ for
some $A\in C(x)$, and every interval $[x,y]$ with $x\in\mathrm{OP}_n$ is Boolean.
\end{proposition}

\begin{proof}
By Proposition~\ref{prop:dictionary}, the elements $y\ge x$ are exactly those with
$C(y)\subseteq C(x)$, and every subset of the chain $C(x)$ is a chain, hence the cut
set of a unique element. So $y\mapsto C(x)\setminus C(y)$ is a bijection from
$[x,\top]$ onto the subsets of $C(x)$. It reverses the inclusion of cut sets, hence
preserves the order in the form $y\le y'\Leftrightarrow
C(x)\setminus C(y)\subseteq C(x)\setminus C(y')$. The statement on covers is
Lemma~\ref{lem:covers}, and an interval $[x,y]$ is an interval of a Boolean lattice.
\end{proof}

Note that the lower intervals $[\bot,y]$ are not Boolean in general: for instance
$\mathcal L_3=[\bot,\top]$ has $14$ elements.

\subsection{Atomicity and coatomicity}

In a lattice with minimum $\bot$ and maximum $\top$, an \emph{atom} is an element
covering $\bot$ and a \emph{coatom} is an element covered by $\top$. The lattice is
\emph{atomic} if every element is a join of atoms, with $\bot$ the empty join, and
\emph{coatomic} if every element is a meet of coatoms, with $\top$ the empty meet. In a
finite lattice every element other than $\top$ lies below some coatom, so this weaker
reading of coatomicity is automatic; what carries content is that every element is a
\emph{meet} of coatoms, and this is the sense used here.

\begin{lemma}[Separating an atom]\label{lem:separate}
Let $x\in\mathrm{OP}_n$ and let $A$ be a proper subset of $[n]$ with $A\notin C(x)$.
Then there is an atom $\rho\le x$ with $A\notin F(\rho)$.
\end{lemma}

\begin{proof}
By Lemma~\ref{lem:atoms}, the atoms below $x=(B_1\mid\cdots\mid B_k)$ are the
concatenations of orderings of $B_1,\dots,B_k$. For such an atom $\rho$, every set
formed by the first $m$ entries of $\rho$ ($0\le m\le n$) equals $A_{j-1}(x)\cup T$ for some $1\le j\le k$, where $T$ is the set
of the first $s$ entries of the ordering of $B_j$ and $0\le s\le|B_j|$; for
$s\in\{0,|B_j|\}$ it is a cut of $x$, or $\emptyset$, or $[n]$.

Suppose first that $A$ is not of the form $A_{j-1}(x)\cup S$ with
$\emptyset\ne S\subsetneq B_j$ for any $j$. Then $A$ is not formed by the first $m$
entries of any atom $\rho\le x$, so $A\notin F(\rho)$ for every such atom, and any one
will do. Otherwise $A=A_{j-1}(x)\cup S$ with $\emptyset\ne S\subsetneq B_j$, and $j$
and $S$ are determined by $A$. Choose the ordering of $B_j$ to begin with an element outside $S$, and any ordering of
the other blocks. In the resulting atom $\rho\le x$, the first $|A|$ entries are
$A_{j-1}(x)$ followed by the first $|S|\ge 1$ entries of the ordering of $B_j$, and
these include an element outside $S$. So the first $|A|$ entries do not form $A$,
that is, $A\notin F(\rho)$.
\end{proof}

\begin{proposition}[Atomic]\label{prop:atomic}
Every element of $\mathcal L_n$ is the join of the atoms below it, and the atoms of
$\mathcal L_n$ are the $n!$ permutations.
\end{proposition}

\begin{proof}
The atoms are the minimal elements of $\mathrm{OP}_n$, which are the permutations by
Lemma~\ref{lem:atoms}. Let $x\in\mathrm{OP}_n$. By Theorem~\ref{thm:lattice}(1) the cut
set of the join of the atoms $\rho\le x$ is the intersection of their flags $F(\rho)$,
which contains $C(x)$ since $F(\rho)\supseteq C(x)$ for every $\rho\le x$. By
Lemma~\ref{lem:separate}, no proper subset $A\notin C(x)$ lies in all of these flags.
So the intersection is $C(x)$, and the join is $x$. Finally $\bot$ is the empty join.
\end{proof}

\begin{corollary}[Coatoms]\label{cor:coatoms}
The coatoms of $\mathcal L_n$ are the elements $(A\mid[n]\setminus A)$ with $A$ a
proper nonempty subset of $[n]$, whose cut set is $\{A\}$. There are $2^n-2$ of
them. The coatoms above $x\in\mathrm{OP}_n$ are exactly those with $A\in C(x)$. A family
$\{A_1\},\dots,\{A_r\}$ of coatoms has meet different from $\bot$ if and only if
$\{A_1,\dots,A_r\}$ is a chain, and then the meet is the element with this cut set.
\end{corollary}

\begin{proof}
An element is a coatom if and only if it has rank $n-2$, that is, exactly one cut.
The coatoms above $x$ are those with $C(\cdot)\subseteq C(x)$. The last statement is
Theorem~\ref{thm:lattice}(2).
\end{proof}

\begin{proposition}[Coatomic]\label{prop:coatomic}
Every element of $\mathcal L_n$ is the meet of the coatoms above it. Namely, if
$x\in\mathrm{OP}_n$ and $x\ne\top$, then $x$ is the meet of the coatoms $\{A\}$ with
$A\in C(x)$; $\top$ is the empty meet; and $\bot$ is the meet of the two coatoms
$\{\{1\}\}$ and $\{\{2\}\}$.
\end{proposition}

\begin{proof}
The union of the cut sets $\{A\}$, $A\in C(x)$, is $C(x)$, which is a chain, so by
Theorem~\ref{thm:lattice}(2) the meet is the element with cut set $C(x)$, which is $x$.
The sets $\{1\}$ and $\{2\}$ are not nested, so their union is not a chain and the
meet is $\bot$.
\end{proof}

\subsection{Pseudo-geometric lattices}

\begin{definition}[Pseudo-geometric lattice]\label{def:pg}
A finite lattice $L$ with minimum $\bot$ and maximum $\top$ is
\emph{pseudo-geometric} if it is graded with $\mathrm{rk}(\bot)=-1$, atomic, and
coatomic in the strong sense above, that is, every element is a meet of coatoms.
\end{definition}

\begin{theorem}\label{thm:pg}
For every $n\ge 2$ the lattice $\mathcal L_n$ is pseudo-geometric. Its atoms are the
$n!$ permutations, it has $2^n-2$ coatoms, and every upper interval $[x,\top]$ with
$x\ne\bot$ is a Boolean lattice.
\end{theorem}

\begin{proof}
Gradedness is Corollary~\ref{cor:graded}, atomicity is
Proposition~\ref{prop:atomic}, coatomicity is Proposition~\ref{prop:coatomic}, and the
upper intervals are Proposition~\ref{prop:boolean}. The counts are those of
Lemma~\ref{lem:atoms} and Corollary~\ref{cor:coatoms}.
\end{proof}
\subsection{Rank modularity and the failure of semimodularity}\label{sec:semimod}

In a geometric lattice the rank is semimodular, $\mathrm{rk}(x\wedge y)+\mathrm{rk}(x\vee y)\le
\mathrm{rk}\,x+\mathrm{rk}\,y$. In $\mathcal L_n$ the rank satisfies the identity with
equality whenever the meet lies above $\bot$, and the defect has no definite sign when the
meet is $\bot$.

\begin{theorem}[Modularity of the rank above $\bot$]\label{thm:modrank}
If $x,y\in\mathcal L_n$ and $x\wedge y>\bot$, then
\[
\mathrm{rk}(x)+\mathrm{rk}(y)=\mathrm{rk}(x\wedge y)+\mathrm{rk}(x\vee y).
\]
\end{theorem}

\begin{proof}
Here $x,y,x\wedge y\in\mathrm{OP}_n$, and $x\vee y\in\mathrm{OP}_n$ by
Theorem~\ref{thm:lattice}. By Theorem~\ref{thm:lattice}(1) and (2),
$C(x\wedge y)=C(x)\cup C(y)$ and $C(x\vee y)=C(x)\cap C(y)$. Since
$\mathrm{rk}=n-1-|C|$ on $\mathrm{OP}_n$, both sides of the claim equal $2(n-1)$ minus
$|C(x)|+|C(y)|$ and minus $|C(x)\cup C(y)|+|C(x)\cap C(y)|$ respectively, and these
agree by inclusion--exclusion.
\end{proof}

\begin{proposition}[The defect when the meet is $\bot$]\label{prop:defect}
Let $x,y\in\mathrm{OP}_n$ with $x\wedge y=\bot$. Then
\[
\mathrm{rk}\,x+\mathrm{rk}\,y-\mathrm{rk}(x\wedge y)-\mathrm{rk}(x\vee y)=n-|C(x)\cup C(y)|.
\]
For every $n\ge3$ both signs occur: two incomparable coatoms give $n-2>0$, and two
opposite atoms $\rho=(1\mid2\mid\cdots\mid n)$ and $\rho'=(n\mid\cdots\mid2\mid1)$ give
$2-n<0$.
\end{proposition}

\begin{proof}
The left-hand side becomes $2(n-1)-|C(x)|-|C(y)|+1-\bigl(n-1-|C(x)\cap C(y)|\bigr)$ by means of $\mathrm{rk}(\bot)=-1$ in the computation of Theorem \ref{thm:modrank}, and this in turn equals $n-|C(x)\cup C(y)|$. Two coatoms $x\ne y$ have $|C(x)|=|C(y)|=1$, and the union of two different singletons has $2$ elements, giving $n-2>0$; the union is not a chain when the two subsets are incomparable, so the meet is $\bot$. The cut sets of $\rho$ and $\rho'$ are the families $\{1,\dots,i\}$ and $\{n-i+1,\dots,n\}$, $1\le i\le n-1$, which are disjoint, so their union has $2(n-1)$ elements and the defect is $2-n$; their union is not a chain for $n\ge3$, so the meet is $\bot$.
\end{proof}

Thus there is no one-sided inequality that could replace semimodularity: the rank is exactly modular where the meet is above $\bot$ and fails in both directions elsewhere. The failure is visible already in the simplest operation, joining an element with an atom.

\begin{proposition}[Joining with an atom]\label{prop:atomjoin}
Let $x\in\mathrm{OP}_n$ and let $\rho$ be an atom with flag $F=C(\rho)$. If $\rho\not\le x$, then $x\wedge\rho=\bot$ and
\[
\mathrm{rk}(x\vee\rho)=\mathrm{rk}(x)+|C(x)\setminus F|.
\]
In particular $x\vee\rho$ covers $x$ if and only if exactly one cut of $x$ is not among the proper prefixes of $\rho$, and the largest possible jump is $n-1-\mathrm{rk}(x)$, in which case $x\vee\rho=\top$. For $n\ge3$ and $\mathrm{rk}(x)\le n-3$ there is an atom $\rho$ such that $x\vee\rho$ does not cover $x$.
\end{proposition}

\begin{proof}
If $C(x)\cup F$ were a chain it would equal $F$, because $F$ is a maximal chain of proper nonempty subsets; then $C(x)\subseteq F$ and $\rho\le x$. So $C(x)\cup F$ is not a chain and $x\wedge\rho=\bot$ by Theorem \ref{thm:lattice}(2). Also,  $\mathrm{rk}(x\vee\rho)=n-1-|C(x)\cap F|=\mathrm{rk}(x)+|C(x)\setminus F|$ by Theorem \ref{thm:lattice}(1). For the last statement, $\mathrm{rk}(x)\le n-3$ means $|C(x)|\ge2$; pick two cuts $A,A'$ of $x$, with $|A|=a$ and $|A'|=a'$. The number of permutations whose prefix of length $a$ is $A$ is $a!(n-a)!$, a fraction $1/\binom na$ of all permutations, and similarly for $A'$. Since
$n\ge3$ and $1\le a,a'\le n-1$, both binomial coefficients are at least $3$, so some permutation $\rho$ has neither $A$ nor $A'$ among its prefixes. Then $\rho\not\le x$ and
$|C(x)\setminus F(\rho)|\ge2$, so the join with $\rho$ raises the rank by at least $2$.
\end{proof}

\begin{example}\label{ex:nocover}
Let $n=4$, $x=(12\mid3\mid4)$ of rank $1$, with $C(x)=\{12,123\}$, and
$\rho=(3\mid1\mid4\mid2)$, whose flag is $\{3,13,134\}$. The two sets are disjoint, so
$x\vee\rho=\top$ has rank $3$, a jump of $2$. For $n=5$, $x=(1\mid2\mid345)$ has rank $2$ and
$C(x)=\{1,12\}$; with $\rho=(2\mid3\mid1\mid4\mid5)$, whose prefixes are
$\{2\},\{2,3\},\{1,2,3\},\{1,2,3,4\}$, the join is $\top$ of rank $4$.
\end{example}

Since $|C(x)|\le n-3$ is needed for a jump of at least $2$ when $\mathrm{rk}(x)\ge2$, a
jump of two or more above an element of rank $\ge2$ first occurs for $n=5$. In this sense
$\mathcal L_n$ is not semimodular: an atom not below $x$ always has meet $\bot$ with it,
and $\bot$ is covered by the atom, but the join does not always cover $x$.

\begin{proposition}[Elements that are modular for the rank]\label{prop:rankmodelts}
Let $n\ge3$ and $x\in\mathcal L_n$. Then $\mathrm{rk}\,x+\mathrm{rk}\,y=\mathrm{rk}(x\wedge y)+
\mathrm{rk}(x\vee y)$ for all $y\in\mathcal L_n$ if and only if $x=\bot$ or $x=\top$.
\end{proposition}

\begin{proof}
For $x=\bot$ and $x=\top$ the identity is trivial. Let $x\in\mathrm{OP}_n$, $x\ne\top$, so that
$C(x)\ne\emptyset$. Choose a cut $B\in C(x)$ and a subset $A\ne B$ of $[n]$ with $|A|=|B|$; it is
incomparable with $B$ and does not belong to $C(x)$. Let $y$ be the coatom with $C(y)=\{A\}$.
The set $C(x)\cup\{A\}$ is not a chain, so $x\wedge y=\bot$, and by
Proposition~\ref{prop:defect} the defect is $n-|C(x)|-1=\mathrm{rk}\,x$. This is nonzero if
$\mathrm{rk}\,x\ge1$. If $x$ is an atom, let $y$ be the atom whose word is the reverse of that of $x$. The
prefixes of $x$ are disjoint from the prefixes of $y$, which are the suffixes of $x$, so
$x\wedge y=\bot$ and the defect is $n-2(n-1)=2-n\ne0$.
\end{proof}

\section{Pseudo-atoms}\label{sec:pseudo}

The atoms and the coatoms of $\mathcal L_n$ were described in Lemma~\ref{lem:atoms}
and Corollary~\ref{cor:coatoms}. In this section we study the elements of rank $1$,
the \emph{pseudo-atoms}, which form the first level above the atoms, and we show
that they, and already the atoms, generate $\mathcal L_n$ with very few elements.
We keep the cut model of Section~\ref{sec:cut}.

\subsection{Pseudo-atoms}

\begin{lemma}[Pseudo-atoms]\label{lem:pseudo}
\begin{enumerate}
\item The pseudo-atoms are the elements of $\mathrm{OP}_n$ with $n-1$ blocks, that
  is, those with one block of size $2$ and all other blocks of size $1$. Each of them
  covers exactly two atoms, which differ by the exchange of two consecutive entries,
  and each atom is covered by exactly $n-1$ pseudo-atoms.
\item There are $\tfrac12(n-1)\,n!$ pseudo-atoms.
\item For $x\in\mathrm{OP}_n$ the number of pseudo-atoms $p\le x$ is
  $\tfrac12\,\mathrm{rk}(x)\prod_{j=1}^k|B_j|!$, where $x=(B_1\mid\cdots\mid B_k)$.
\end{enumerate}
\end{lemma}

\begin{proof}
(1) The first sentence is the definition $\mathrm{rk}(p)=n-k(p)=1$. By
Lemma~\ref{lem:atoms}, the atoms below $p$ number $\prod_j|B_j|!=2$, and they are the
two orderings of the block $\{a,b\}$ of size $2$, the other blocks being singletons;
they differ by the exchange of the two consecutive entries $a,b$. An atom $\rho$ has
$C(\rho)=F(\rho)$ with $n-1$ cuts, and by Lemma~\ref{lem:covers} the elements covering
$\rho$ are obtained by deleting one cut, so there are $n-1$ of them.

(2) Count the pairs (atom, pseudo-atom covering it) in two ways: there are $n!(n-1)$ of
them by the second statement of (1) and twice as many as there are pseudo-atoms by the
first.

(3) Count the pairs $(\rho,p)$ with $\rho$ an atom, $p$ a pseudo-atom and
$\rho<p\le x$. By Lemma~\ref{lem:atoms} there are $\prod_j|B_j|!$ atoms $\rho\le x$.
For each of them the interval $[\rho,x]$ is Boolean on $F(\rho)\setminus C(x)$
(Proposition~\ref{prop:boolean}), a set with $n-1-|C(x)|=\mathrm{rk}(x)$ elements, so
$\rho$ is covered in $[\rho,x]$ by $\mathrm{rk}(x)$ pseudo-atoms. On the other hand every
pseudo-atom $p\le x$ covers exactly two atoms, and both lie below $x$. Hence
$2\cdot\#\{p\le x\}=\mathrm{rk}(x)\prod_j|B_j|!$.
\end{proof}

For $x=\top$ the last count gives $\tfrac12(n-1)\,n!$ again. For $n=3$ the
pseudo-atoms are the six elements of rank $1$ in Example~\ref{ex:n3}; they coincide
with the coatoms, as $n-2=1$.

\subsection{Generating an element by two atoms}

For $x=(B_1\mid\cdots\mid B_k)\in\mathrm{OP}_n$ let $\rho_x^{+}$ and $\rho_x^{-}$ be
the atoms obtained by listing the blocks $B_1,\dots,B_k$ in this order, each block
written in increasing order for $\rho_x^{+}$ and in decreasing order for
$\rho_x^{-}$. Both lie below $x$ by Lemma~\ref{lem:atoms}. For instance, if
$x=(12\mid 34)$ then $\rho_x^{+}=(1\mid2\mid3\mid4)$ and
$\rho_x^{-}=(2\mid1\mid4\mid3)$.

\begin{proposition}[Two atoms suffice]\label{prop:twoatoms}
Let $x\in\mathrm{OP}_n$ with $\mathrm{rk}(x)\ge1$. Then $\rho_x^{+}\ne\rho_x^{-}$ and
$x=\rho_x^{+}\vee\rho_x^{-}$. In particular every element of $\mathrm{OP}_n$ of rank at
least $1$ is a join of two distinct atoms, while an atom is not.
\end{proposition}

\begin{proof}
Since $\mathrm{rk}(x)\ge1$ some block has size at least $2$, and its two orderings
differ, so $\rho_x^{+}\ne\rho_x^{-}$. By Theorem~\ref{thm:lattice}(1) we have
$C(\rho_x^{+}\vee\rho_x^{-})=F(\rho_x^{+})\cap F(\rho_x^{-})$, and this contains
$C(x)$ because both atoms lie below $x$. Conversely, let $1\le m\le n-1$ and write
$m=|A_{j-1}(x)|+s$ with $1\le j\le k$ and $0\le s<|B_j|$. The set of the first $m$
entries of $\rho_x^{+}$ is $A_{j-1}(x)\cup S^{+}$, where $S^{+}$ consists of the $s$
smallest elements of $B_j$, and that of $\rho_x^{-}$ is $A_{j-1}(x)\cup S^{-}$, where
$S^{-}$ consists of the $s$ largest. If $0<s<|B_j|$ then the smallest element of $B_j$
lies in $S^{+}$ but not in $S^{-}$, so the two sets differ and neither is in the
intersection of the flags. If $s=0$ the common set is $A_{j-1}(x)$, which is a cut of
$x$ (note $j\ge2$ since $m\ge1$). Hence
$F(\rho_x^{+})\cap F(\rho_x^{-})=C(x)$ and the join is $x$.

For the last statement, if $\rho\ne\rho'$ are atoms then $\rho'\not\le\rho$, so
$\rho\vee\rho'\ne\rho$ and the join is not an atom.
\end{proof}

\begin{corollary}[Unique join representations]\label{cor:unique}
For $x\in\mathcal L_n$ the following are equivalent: $\mathrm{rk}(x)\le1$, and there
is exactly one set of atoms whose join is $x$. Thus, apart from $\bot$ and the atoms,
the pseudo-atoms are the only elements with a unique representation as a join of atoms.
\end{corollary}

\begin{proof}
By Proposition~\ref{prop:atomic}, the set of all atoms below $x$ has join $x$. For
$x=\bot$ no nonempty set of atoms has join $\bot$, so the empty set is the only
representation. If $x=\rho$ is an atom, every atom in a set with join $\rho$ lies below
$\rho$, hence equals $\rho$, and the only representation is $\{\rho\}$. If $x=p$ is a
pseudo-atom, exactly two atoms $\rho,\rho'$ lie below $p$ (Lemma~\ref{lem:pseudo}), and
the proper subsets of $\{\rho,\rho'\}$ have joins $\bot$, $\rho$, $\rho'$, none equal to
$p$; so $\{\rho,\rho'\}$ is the only representation. Finally let $\mathrm{rk}(x)\ge2$,
with $x=(B_1\mid\cdots\mid B_k)$. Since $\mathrm{rk}(x)=\sum_j(|B_j|-1)$, either some
block has at least $3$ elements or at least two blocks have at least $2$, so
$\prod_j|B_j|!\ge4$ atoms lie below $x$. The set of all of them and the set
$\{\rho_x^{+},\rho_x^{-}\}$ of Proposition~\ref{prop:twoatoms} are two different
representations of $x$.
\end{proof}

The contrast with geometric lattices is clear. In a geometric lattice an element of
rank $r$ is a join of $r$ atoms and of no fewer; here two atoms suffice whatever the
rank, which is at most $n-1$. For example, in $\mathcal L_4$ the atoms
$(1\mid2\mid3\mid4)$ and $(4\mid3\mid2\mid1)$ have disjoint flags, so their join is
$\top$, of rank $3$.

\subsection{Generating an element by two pseudo-atoms}

\begin{proposition}[Two pseudo-atoms suffice]\label{prop:twopseudo}
Let $x\in\mathrm{OP}_n$ with $\mathrm{rk}(x)\ge1$. There are pseudo-atoms $p,q\le x$
with $x=p\vee q$, and one can take $p\ge\rho_x^{+}$ and $q\ge\rho_x^{-}$. If
$\mathrm{rk}(x)\ge2$ then no single pseudo-atom suffices, so two is the minimum number
of pseudo-atoms whose join is $x$.
\end{proposition}

\begin{proof}
Since $\rho_x^{+}<x$ the interval $[\rho_x^{+},x]$ has rank $\mathrm{rk}(x)\ge1$, so
it contains an element $p$ covering $\rho_x^{+}$; by Lemma~\ref{lem:covers}, $p$ has
rank $1$. Choose $q$ in the same way above $\rho_x^{-}$. Then
$C(x)\subseteq C(p)\cap C(q)\subseteq F(\rho_x^{+})\cap F(\rho_x^{-})=C(x)$, where the
first inclusion holds since $p,q\le x$, the second since $p\ge\rho_x^{+}$ and
$q\ge\rho_x^{-}$, and the equality is Proposition~\ref{prop:twoatoms}. Thus
$C(p\vee q)=C(x)$ and $p\vee q=x$. If $\mathrm{rk}(x)\ge2$, a single pseudo-atom has
rank $1$ and is not $x$.
\end{proof}

\begin{example}\label{ex:twopseudo}
For $x=(12\mid34)$ and the atoms above, take the pseudo-atoms obtained by deleting the first cut of each atom,
$p=(12\mid3\mid4)\ge(1\mid2\mid3\mid4)$
and $q=(12\mid4\mid3)\ge(2\mid1\mid4\mid3)$. Their cut sets are $\{12,123\}$ and
$\{12,124\}$, with intersection $\{12\}=C(x)$, so $p\vee q=x$.
\end{example}

\begin{corollary}[Pseudo-atoms generate]\label{cor:pseudogen}
Every element $x$ of $\mathcal L_n$ with $\mathrm{rk}(x)\ge1$ is the join of the
pseudo-atoms below it; in fact two of them already suffice.
\end{corollary}

\begin{proof}
By Proposition~\ref{prop:twopseudo}, $x=p\vee q$ with $p,q\le x$ pseudo-atoms. The join
of all the pseudo-atoms below $x$ lies between $p\vee q$ and $x$.
\end{proof}

The elements of rank $0$ are atoms and are not joins of pseudo-atoms, and $\bot$ is the
empty join; so the pseudo-atoms generate the part of $\mathcal L_n$ of rank at least
$1$, in the same way as the atoms generate all of $\mathcal L_n$ (Proposition~\ref{prop:atomic}).

\begin{remark}[Redundancy of join representations]\label{rem:redundancy}
By Proposition~\ref{prop:atomic}, $x$ is the join of all the atoms below it, and there
are $\prod_j|B_j|!$ of them; this representation is highly redundant. Call a family of
atoms with join $x$ \emph{irredundant} if no proper subfamily has join $x$. By
Proposition~\ref{prop:twoatoms} two atoms suffice when $\mathrm{rk}(x)\ge1$, but
a representation is unique only when $\mathrm{rk}(x)\le1$ (Corollary~\ref{cor:unique}), and irredundant families are in general neither unique nor of constant size. For $n=3$ we have
$\top=(1\mid2\mid3)\vee(2\mid3\mid1)=(1\mid2\mid3)\vee(3\mid2\mid1)$, since the flag
$\{1,12\}$ of the first atom is disjoint from the flags $\{2,23\}$ and $\{3,23\}$ of the
other two. For
$n=4$ the atoms $(1\mid2\mid3\mid4)$, $(1\mid3\mid4\mid2)$, $(3\mid1\mid2\mid4)$ have
flags $\{1,12,123\}$, $\{1,13,134\}$, $\{3,13,123\}$, whose pairwise intersections are
nonempty, $\{1\}$, $\{123\}$, $\{13\}$, while the intersection of all three is empty. So
$\top$ is the join of an irredundant family of three atoms, although two suffice. In a
geometric lattice, by contrast, all irredundant families of atoms with a given join
have the same size, the rank of the join, by the basis exchange property of
matroids. However, the maximal (redundant) join representation of an element, the set of all atoms below it, is unique, and its size is the lattice function whose Möbius inversion takes value $1$ on atoms and $0$ elsewhere, as in Rossi (2018) \cite{Rossi18}. 
We write $s(x)=|\mathrm{At}(x)|=\prod_j|B_j|!$ for the size of the maximal
representation, as in Section~\ref{sec:metric}. Let $m(x)$ be the least size of a
family of atoms with join $x$, so that $m(\bot)=0$, $m(\rho)=1$ for an atom $\rho$, and
$m(x)=2$ if $\mathrm{rk}(x)\ge1$, by Proposition~\ref{prop:twoatoms}. The
\emph{redundancy} of $x$ is the number of atoms that can be discarded from the maximal
representation, $\mathrm{red}(x)=s(x)-m(x)$. Thus $\mathrm{red}(\top)=n!-2$,
$\mathrm{red}(x)=\prod_j|B_j|!-2$ if $\mathrm{rk}(x)\ge1$, and by
Corollary~\ref{cor:unique} $\mathrm{red}(x)=0$ if and only if $\mathrm{rk}(x)\le1$.
Irredundant families can be larger than the minimum, as in the example with three atoms
above, and the largest size $c(x)$ of an irredundant family is $\mathrm{rk}(x)$ when
$\mathrm{rk}(x)\ge2$ (Proposition~\ref{prop:irredundant}).
\end{remark}

\begin{remark}[Geometry; not used later]
By Remark~\ref{rem:geometry}, the pseudo-atoms are the
edges of the permutohedron, and Lemma~\ref{lem:pseudo}(3) is the count of edges of a
face with $N$ vertices and dimension $d$, namely $Nd/2$, since the graph of a face is
$d$-regular. Again we do not use this identification.
\end{remark}

\section{Size and distances}\label{sec:metric}

In this section we attach two kinds of distances to $\mathcal L_n$. The first uses the
atoms below an element and is the Hamming distance between sets of atoms; it is
defined on the whole lattice. The second uses the cuts, hence the coatoms above an
element; on atoms it is the rank of the join, so that it measures how far apart two atoms
are by how many cuts they do not share. Both are metrics that embed
isometrically in a hypercube.

\subsection{The size of an element}

For $x\in\mathcal L_n$ let $\mathrm{At}(x)$ be the set of atoms $\rho\le x$ and put
\[
s(x)=|\mathrm{At}(x)|.
\]
By Lemma~\ref{lem:atoms}, $s(x)=\prod_{j=1}^k|B_j|!$ for $x=(B_1\mid\cdots\mid B_k)
\in\mathrm{OP}_n$, so $s(\rho)=1$ for an atom, $s(\top)=n!$, and $s(\bot)=0$. By
Proposition~\ref{prop:atomic}, $s(x)$ is the size of the maximal join representation
of $x$ in Remark~\ref{rem:redundancy}, which is the family $\mathrm{At}(x)$ of all
atoms below $x$.

\begin{proposition}[Atom sets]\label{prop:atomsets}
Let $x,y\in\mathcal L_n$.
\begin{enumerate}
\item $x\le y$ if and only if $\mathrm{At}(x)\subseteq\mathrm{At}(y)$. In particular
  $x\mapsto\mathrm{At}(x)$ is injective.
\item $\mathrm{At}(x\wedge y)=\mathrm{At}(x)\cap\mathrm{At}(y)$ and
  $\mathrm{At}(x\vee y)\supseteq\mathrm{At}(x)\cup\mathrm{At}(y)$.
\item $s(x)=\sum_{y\le x}f(y)$, where $f(y)=1$ if $y$ is an atom and $f(y)=0$
  otherwise. Hence the Möbius inversion of $s$ is $f$, and it is nonnegative.
\end{enumerate}
\end{proposition}

\begin{proof}
(1) If $x\le y$ every atom below $x$ is below $y$, while if
$\mathrm{At}(x)\subseteq\mathrm{At}(y)$, then 
$x=\bigvee\mathrm{At}(x)\le\bigvee\mathrm{At}(y)=y$ by Proposition \ref{prop:atomic}. (2) An atom is below $x\wedge y$
if and only if it is below $x$ and below $y$, and $x\vee y$ is above both $x$ and $y$.
(3) The sum counts the atoms $y\le x$, and Möbius inversion on the finite poset
$\mathcal L_n$ is unique.
\end{proof}

The function $s$ is the lattice analogue of the cardinality of a subset, and the
Möbius inversion in (3) is the one used in \cite{Rossi18}. It has more than the
supermodularity one would expect.

\begin{proposition}[Total monotonicity]\label{prop:monotone}
For every $x_1,\dots,x_k\in\mathcal L_n$,
\[
s\Bigl(\bigvee_{i=1}^kx_i\Bigr)\ \ge\ \sum_{\emptyset\ne I\subseteq[k]}(-1)^{|I|+1}
s\Bigl(\bigwedge_{i\in I}x_i\Bigr).
\]
In particular, $s$ is supermodular:
$s(x\wedge y)+s(x\vee y)\ge s(x)+s(y)$, with equality if and only if every atom below
$x\vee y$ lies below $x$ or below $y$.
\end{proposition}

\begin{proof}
By Proposition~\ref{prop:atomsets}(2), the left-hand side is at least
$\bigl|\bigcup_i\mathrm{At}(x_i)\bigr|$ and $s(\bigwedge_{i\in I}x_i)=
\bigl|\bigcap_{i\in I}\mathrm{At}(x_i)\bigr|$. The right-hand side is then the
inclusion--exclusion expression of $\bigl|\bigcup_i\mathrm{At}(x_i)\bigr|$. For $k=2$,
equality holds if and only if $\mathrm{At}(x\vee y)=\mathrm{At}(x)\cup\mathrm{At}(y)$.
\end{proof}

Both statements hold in every finite atomic lattice, with the same proofs.

\begin{corollary}[Lower intervals]\label{cor:boolean}
For every $x\in\mathcal L_n$ one has $s(x)\ge\mathrm{rk}(x)+1$, and the following are
equivalent:
\begin{enumerate}
\item[(a)] $\mathrm{rk}(x)\le1$;
\item[(b)] $x$ has exactly one representation as a join of atoms;
\item[(c)] $s(x)=\mathrm{rk}(x)+1$;
\item[(d)] the lower interval $[\bot,x]$ is a Boolean lattice.
\end{enumerate}
\end{corollary}

\begin{proof}
For $x=\bot$ we have $s(x)=0=\mathrm{rk}(x)+1$. Let $x=(B_1\mid\cdots\mid B_k)$ and
$r=\mathrm{rk}(x)=\sum_j(|B_j|-1)$. Since $b!\ge2^{b-1}$ for $b\ge1$, we get
$s(x)=\prod_j|B_j|!\ge2^r$. Now $2^r\ge r+1$ with equality exactly for $r\in\{0,1\}$,
and $2^r\ge r+2$ for $r\ge2$; hence $s(x)\ge\mathrm{rk}(x)+1$, with equality if and only
if $\mathrm{rk}(x)\le1$. This is (a)$\Leftrightarrow$(c), and (a)$\Leftrightarrow$(b)
is Corollary~\ref{cor:unique}. If (a) holds, $[\bot,x]$ is $\{\bot\}$, or $\{\bot,x\}$
for an atom, or $\{\bot,\rho,\rho',x\}$ for a pseudo-atom, with $\rho,\rho'$ the two
atoms below it (Lemma~\ref{lem:pseudo}); each is Boolean. Finally, if $[\bot,x]$ is
Boolean, its rank is $\mathrm{rk}(x)+1$ and it has that many atoms, which are the atoms
of $\mathcal L_n$ below $x$; so $s(x)=\mathrm{rk}(x)+1$, which is (c).
\end{proof}

In a geometric lattice the number of atoms below $x$ is at least the rank of $x$, with
equality exactly when $[\bot,x]$ is Boolean; here the same holds with the shift by one
due to $\mathrm{rk}(\bot)=-1$.

\subsection{The closure operator on atoms}

For a set $S$ of atoms put $\mathrm{cl}(S)=\mathrm{At}\bigl(\bigvee S\bigr)$, the set of
atoms below the join of $S$, with $\mathrm{cl}(\emptyset)=\mathrm{At}(\bot)=\emptyset$.

\begin{proposition}[Closure operator]\label{prop:closure}
\begin{enumerate}
\item $\mathrm{cl}$ is a closure operator on the subsets of the set of atoms:
  $S\subseteq\mathrm{cl}(S)$; $S\subseteq T$ implies $\mathrm{cl}(S)\subseteq\mathrm{cl}(T)$;
  $\mathrm{cl}(\mathrm{cl}(S))=\mathrm{cl}(S)$.
\item The closed sets are exactly the sets $\mathrm{At}(x)$, $x\in\mathcal L_n$, and
  $x\mapsto\mathrm{At}(x)$ is an isomorphism of $\mathcal L_n$ onto the lattice of closed
  sets ordered by inclusion, in which the meet is the intersection and the join of two
  closed sets is the closure of their union.
\item $\mathrm{cl}$ does not satisfy the exchange property.
\end{enumerate}
\end{proposition}

\begin{proof}
(1) The first two properties are clear. By Proposition~\ref{prop:atomic},
$\bigvee\mathrm{cl}(S)=\bigvee\mathrm{At}(\bigvee S)=\bigvee S$, which gives the third.
(2) By definition $\mathrm{cl}(S)$ is of the form $\mathrm{At}(x)$, and
$\mathrm{cl}(\mathrm{At}(x))=\mathrm{At}(x)$ because $\bigvee\mathrm{At}(x)=x$. The map
is an order isomorphism onto its image by Proposition~\ref{prop:atomsets}(1); moreover
$\mathrm{At}(x\wedge y)=\mathrm{At}(x)\cap\mathrm{At}(y)$ and
$\mathrm{At}(x\vee y)=\mathrm{cl}(\mathrm{At}(x)\cup\mathrm{At}(y))$, since
$\bigvee(\mathrm{At}(x)\cup\mathrm{At}(y))=x\vee y$.
(3) Take $n=3$, $S=\{(1\mid2\mid3)\}$, $p=(1\mid3\mid2)$ and $q=(3\mid2\mid1)$. The flags
of $(1\mid2\mid3)$ and $q$ are disjoint, so their join is $\top$ and
$p\in\mathrm{cl}(S\cup\{q\})$, while $p\notin\mathrm{cl}(S)=S$. On the other hand the
flags $\{1,12\}$ and $\{1,13\}$ of $(1\mid2\mid3)$ and $p$ meet in $\{1\}$, so their join
is $(1\mid23)$, which lies above only these two atoms; hence $q\notin
\mathrm{cl}(S\cup\{p\})$.
\end{proof}

A finite lattice is geometric if and only if it is the lattice of closed sets of a
closure operator on its atoms which satisfies the exchange property, that is, of a simple
matroid \cite{StanleyEC1,Oxley}. Thus $\mathcal L_n$ is the lattice of closed sets of an atomistic
closure operator, and what separates it from the geometric case is precisely the failure
of exchange in part (3).

\subsection{Irredundant families of atoms}

A family of atoms is \emph{irredundant for $x$} if its join is $x$ and no proper
subfamily has join $x$; in the language of the closure operator, it is a minimal
generating set of the closed set $\mathrm{At}(x)$. As in Remark~\ref{rem:redundancy}, let
$m(x)$ and $c(x)$ be the smallest and the largest size of an irredundant family for $x$;
the smallest size of any family with join $x$ is $m(x)$. In a geometric lattice all
irredundant families of atoms for $x$ have the same size, the rank of $x$. Here they do
not, but their sizes are controlled. (These are the families called independent by Sagan \cite{Sagan95} and in \cite[Section~9]{BlassSagan}: the join of every proper subfamily is smaller than the join of the family. They are not the independent sets of Raz of Section~\ref{sec:open}.)

\begin{proposition}[Sizes of irredundant families]\label{prop:irredundant}
For $x\in\mathcal L_n$ one has $m(x)=c(x)=0$ if $x=\bot$, $m(x)=c(x)=1$ if $x$ is an
atom, and $m(x)=c(x)=2$ if $\mathrm{rk}(x)=1$. If $\mathrm{rk}(x)\ge2$, then
$m(x)=2$ and $c(x)=\mathrm{rk}(x)$.
\end{proposition}

\begin{proof}
For $\mathrm{rk}(x)\le1$ the family of atoms with join $x$ is unique, by
Corollary~\ref{cor:unique}, and the values follow. Let $\mathrm{rk}(x)\ge2$. Then
$m(x)=2$ by Proposition~\ref{prop:twoatoms}, since one atom has rank $0$. Write $r=
\mathrm{rk}(x)$, and for an atom $\rho\le x$ write $F(\rho)\supseteq C(x)$ for its flag.

\emph{Upper bound.} Let $\rho_1,\dots,\rho_k$ be an irredundant family for $x$ with
$k\ge3$ (a family with $k\le2$ has $k\le r$), and put $F_i=F(\rho_i)$, so that
$\bigcap_iF_i=C(x)$ by Theorem~\ref{thm:lattice}(1). Irredundancy means that
$\bigcap_{j\ne i}F_j\supsetneq C(x)$ for every $i$, so we can choose a cut
$a_i\in\bigcap_{j\ne i}F_j\setminus C(x)$. Then $a_i\notin F_i$, since otherwise $a_i$
would lie in $\bigcap_jF_j=C(x)$; hence the $a_i$ are distinct, and $a_j\in F_i$ for
$j\ne i$. Any two of them, $a_p$ and $a_q$, lie in a common flag $F_i$ with $i\ne p,q$,
which exists because $k\ge3$; so they are comparable. Each $a_j$ also lies in some
$F_i\supseteq C(x)$, so it is comparable with every element of $C(x)$. Therefore
$C(x)\cup\{a_1,\dots,a_k\}$ is a chain of proper subsets of $[n]$ with $|C(x)|+k$
elements, and $|C(x)|+k\le n-1$, that is, $k\le r$.

\emph{Lower bound.} Let $\pi$ be an atom below $x$, for instance the one listing the
blocks of $x$ in order, each in increasing order, and let $P_t$ be the set of the
first $t$ entries of $\pi$. The $t$ with $1\le t\le n-1$ and $P_t\notin C(x)$ are exactly
$r$ in number, say $t\in T$, and for such a $t$ the entries $\pi_t,\pi_{t+1}$ lie in the
same block of $x$. Let $\rho_t$ be $\pi$ with these two entries exchanged. By
Lemma~\ref{lem:atoms}, $\rho_t\le x$, and its flag is $F(\pi)$ with $P_t$ replaced by
the set $P_t'=(P_t\setminus\{\pi_t\})\cup\{\pi_{t+1}\}\ne P_t$. The sets $P_t$ and $P_t'$
are the only ones of size $t$ that occur in the flags $F(\pi)$ and $F(\rho_t)$. Now
$\bigcap_{t\in T}F(\rho_t)=C(x)$: a cut of $x$ lies in every flag; $P_s$ with $s\in T$
is not in $F(\rho_s)$; and $P_s'$ is not in $F(\rho_t)$ for $t\ne s$, and such a $t$
exists because $r\ge2$. So the join of the family $\{\rho_t:t\in T\}$ is $x$. If
$\rho_s$ is removed, the intersection of the remaining flags contains $P_s\notin
C(x)$, so the join is not $x$. The family is irredundant with $r$ elements, and
$c(x)\ge r$.
\end{proof}

For example, the $n-1$ atoms obtained from $(1\mid2\mid\cdots\mid n)$ by exchanging two
consecutive entries form an irredundant family with join $\top$, and $c(\top)=n-1$. For
$n=4$ and $x=\top$ the sizes $2$ and $3$ both occur, as in
Remark~\ref{rem:redundancy}. The contrast with geometric lattices is that there the size
of an irredundant family is the rank of $x$ counted from the atoms, which here is
$\mathrm{rk}(x)+1$; the largest irredundant family is one atom short of it, and the
smallest has only two.

\subsection{The Hamming distance between atom sets}

\begin{definition}\label{def:d}
For $x,y\in\mathcal L_n$ put
$d(x,y)=s(x)+s(y)-2\,s(x\wedge y)$.
\end{definition}

For subsets $S,T$ of a finite set, $|S\triangle T|=|S|+|T|-2|S\cap T|$ is the
Hamming distance between their indicator vectors. By
Proposition~\ref{prop:atomsets}(2), $d$ is exactly this distance for the atom sets.

\begin{theorem}[Hamming metric]\label{thm:hamming}
\begin{enumerate}
\item $d(x,y)=|\mathrm{At}(x)\triangle\mathrm{At}(y)|$, and $d$ is a metric on
  $\mathcal L_n$. The map $x\mapsto\mathrm{At}(x)$ is an isometric embedding of
  $(\mathcal L_n,d)$ in the Hamming cube $\{0,1\}^{n!}$ indexed by the atoms.
\item Give the edge $u\lessdot v$ of the Hasse diagram of $\mathcal L_n$ the weight
  $s(v)-s(u)$. Then $d(x,y)$ is the minimum weight of a path from $x$ to $y$, and it
  is attained by a path going down from $x$ to $x\wedge y$ and then up to $y$.
\end{enumerate}
\end{theorem}

\begin{proof}
(1) The formula follows from $s(x\wedge y)=|\mathrm{At}(x)\cap\mathrm{At}(y)|$. The
triangle inequality is that of the Hamming distance, and $d(x,y)=0$ forces
$\mathrm{At}(x)=\mathrm{At}(y)$, hence $x=y$ by Proposition~\ref{prop:atomsets}(1).

(2) For $u\lessdot v$ we have $\mathrm{At}(u)\subseteq\mathrm{At}(v)$, so the weight
$s(v)-s(u)$ equals $d(u,v)$. By the triangle inequality the weight of any path from
$x$ to $y$ is at least $d(x,y)$. A saturated chain from $x$ down to $m=x\wedge y$ has
weight $s(x)-s(m)$, since the weights telescope, and a saturated chain from $m$ up to
$y$ has weight $s(y)-s(m)$. Their sum is $d(x,y)$.
\end{proof}

By supermodularity, since
$\bigl(2s(x\vee y)-s(x)-s(y)\bigr)-\bigl(s(x)+s(y)-2s(x\wedge y)\bigr)\ge0$, a path through the meet is never longer than one through the join.

For $x,y\in\mathrm{OP}_n$ with $C(x)\cup C(y)$ a chain,
\[
d(x,y)=\prod_j|B_j(x)|!+\prod_j|B_j(y)|!-2\prod_j|B_j(x\wedge y)|!,
\]
and otherwise $d(x,y)=s(x)+s(y)$.

\begin{example}\label{ex:d}
Two distinct atoms have $d=2$, since their meet is $\bot$; so $d$ does not
distinguish between pairs of atoms. Two pseudo-atoms have $d=2$ if they cover a
common atom and $d=4$ otherwise: for $n=4$ the pseudo-atoms $(12\mid3\mid4)$ and
$(13\mid2\mid4)$ have cut sets $\{12,123\}$ and $\{13,123\}$, whose union is not a
chain, so their meet is $\bot$ and $d=2+2=4$. Finally $d(\rho,\top)=n!-1$ for every
atom $\rho$.
\end{example}

\subsection{Cuts and the distance between atoms}

For atoms the metric $d$ is the discrete metric, so we need another one.

\begin{definition}\label{def:delta}
For atoms $\rho,\rho'$ put $\delta(\rho,\rho')=\mathrm{rk}(\rho\vee\rho')$.
\end{definition}

By Definition~\ref{def:rank} and Theorem~\ref{thm:lattice}(1),
$\delta(\rho,\rho')=n-1-|F(\rho)\cap F(\rho')|$, since $C(\rho\vee\rho')=
F(\rho)\cap F(\rho')$. Here $F(\rho)\cap F(\rho')$ is the set of proper subsets which
are the first $m$ entries of both $\rho$ and $\rho'$ for some $m$, and $\delta$ counts
the values $m$, $1\le m\le n-1$, for which this does not happen.

\begin{theorem}[The metric $\delta$]\label{thm:delta}
Let $\rho,\rho',\rho''$ be atoms.
\begin{enumerate}
\item $\delta(\rho,\rho')=\tfrac12\,|F(\rho)\triangle F(\rho')|$. Hence $\delta$ is a
  metric on the set of atoms, and $\rho\mapsto F(\rho)$ is an isometric embedding of
  $(\{\text{atoms}\},2\delta)$ in the Hamming cube $\{0,1\}^{2^n-2}$ indexed by the
  coatoms.
\item $0\le\delta\le n-1$, and $\delta(\rho,\rho')=n-1$ if and only if
  $F(\rho)\cap F(\rho')=\emptyset$, that is, $\rho\vee\rho'=\top$. This happens for
  $\rho$ and its reversal.
\item $\delta(\rho,\rho')=1$ if and only if $\rho,\rho'$ differ by the exchange of two
  consecutive entries, that is, if they are the two atoms below a pseudo-atom.
\end{enumerate}
\end{theorem}

\begin{proof}
(1) Both flags have $n-1$ elements, so
$|F(\rho)\triangle F(\rho')|=2(n-1)-2|F(\rho)\cap F(\rho')|=2\delta(\rho,\rho')$. The
rest is the Hamming distance. (2) is immediate from the formula. (3) Two distinct
atoms below a pseudo-atom have $\delta=1$ (Lemma~\ref{lem:pseudo}); conversely
$\delta=1$ means $\rho\vee\rho'$ has rank $1$, hence is a pseudo-atom with $\rho,\rho'$
below it.
\end{proof}

There is a coatom-side analogue of $d$ for all elements of $\mathrm{OP}_n$, and it is
the graph distance in the Hasse diagram.

\begin{proposition}[Cut distance]\label{prop:dstar}
For $x,y\in\mathrm{OP}_n$ put $d^{*}(x,y)=|C(x)\triangle C(y)|$. Then
\[
d^{*}(x,y)=2\,\mathrm{rk}(x\vee y)-\mathrm{rk}(x)-\mathrm{rk}(y),
\]
$d^{*}$ is a metric on $\mathrm{OP}_n$, and $d^{*}(x,y)$ is the distance between $x$
and $y$ in the Hasse diagram of $\mathrm{OP}_n$ with all edges of length $1$. On atoms,
$d^{*}=2\delta$.
\end{proposition}

\begin{proof}
We have $|C(x)\triangle C(y)|=|C(x)|+|C(y)|-2|C(x)\cap C(y)|$ with
$|C(x)|=n-1-\mathrm{rk}(x)$ and $C(x)\cap C(y)=C(x\vee y)$, which gives the formula. The
map $x\mapsto C(x)$ is injective (Proposition~\ref{prop:dictionary}), so $d^{*}$ is a
metric. Covers change the cut set by one element, so a path from $x$ to $y$ has length
at least $d^{*}(x,y)$. The upper intervals are Boolean (Proposition~\ref{prop:boolean}),
so there is a saturated chain from $x$ up to $x\vee y$ of length
$\mathrm{rk}(x\vee y)-\mathrm{rk}(x)$, and similarly for $y$; the path through $x\vee
y$ has length $d^{*}(x,y)$.
\end{proof}

Thus $d$ and $d^{*}$ are the two sides of the same picture: the Hamming distance for
the atoms below an element, which is weighted by size and goes through the meet, and
the Hamming distance for the coatoms above an element, which is the unweighted
graph distance and goes through the join.

\subsection{The inversion graph}

Let $\rho,\rho'$ be atoms, regarded as orderings of $[n]$. Say that $a,b\in[n]$ are
\emph{inverted} if they appear in different relative order in $\rho$ and $\rho'$, and
let $\mathcal G(\rho,\rho')$ be the graph on $[n]$ whose edges are the inverted pairs.
The number of edges is the Kendall distance $\tau(\rho,\rho')$, the number of
adjacent exchanges needed to go from $\rho$ to $\rho'$. If $\rho$ is the identity,
$\mathcal G$ is the inversion graph of the permutation $\rho'$.

\begin{theorem}[Join of two atoms]\label{thm:components}
The blocks of $\rho\vee\rho'$ are the vertex sets of the connected components of
$\mathcal G(\rho,\rho')$. Consequently
\[
\delta(\rho,\rho')=n-c\bigl(\mathcal G(\rho,\rho')\bigr)\ \le\ \tau(\rho,\rho'),
\]
where $c$ is the number of connected components, with equality if and only if
$\mathcal G(\rho,\rho')$ is a forest.
\end{theorem}

\begin{proof}
Let $P_0=\emptyset\subset P_1\subset\cdots\subset P_l=[n]$ be the common prefix sets
of $\rho$ and $\rho'$, together with $\emptyset$ and $[n]$, that is, the cuts of
$\rho\vee\rho'$ with the two ends. The blocks of $\rho\vee\rho'$ are the sets
$D_j=P_j\setminus P_{j-1}$.

\emph{Each edge lies inside a block.} Let $a,b$ be inverted, with $a$ before $b$ in
$\rho$ and $b$ before $a$ in $\rho'$. Suppose a common prefix set $P$ contained exactly
one of them. As a prefix of $\rho$ it contains $a$ and not $b$, and as a prefix of
$\rho'$ it contains $b$ if it contains $a$; this is a contradiction. So no $P_j$
separates $a$ and $b$.

\emph{Each block is connected.} Fix $j$ and suppose that the subgraph on $D=D_j$ is not
connected. Both $\rho$ and $\rho'$ list the elements of $D$ consecutively, after
$P_{j-1}$. Let $K$ be the component of $\mathcal G$ on $D$ containing the element $u$ of
$D$ that comes first in $\rho$, and let $c\in D\setminus K$. Then $c$ is not inverted
with any element of $K$, so it has the same relative position with respect to every
element of $K$ in $\rho$ and in $\rho'$. We claim that every $k\in K$ is before $c$ in
$\rho$. This holds for $u$. If it failed for some $k$, take a path in $K$ from $u$ to
$k$; then some consecutive pair $k_i,k_{i+1}$ on it has $k_i$ before $c$ and $c$ before
$k_{i+1}$ in $\rho$, hence in $\rho'$ too, so $k_i$ is before $k_{i+1}$ in both
orderings, contradicting that they are inverted. Therefore in both $\rho$ and $\rho'$ the
elements of $K$ precede the elements of $D\setminus K$, and $P_{j-1}\cup K$ is a common
prefix set strictly between $P_{j-1}$ and $P_j$, which is impossible.

Hence blocks and components coincide, and $\delta=\mathrm{rk}(\rho\vee\rho')=n-k$ with
$k$ the number of blocks. Finally a graph on $n$ vertices with $c$ components has at
least $n-c$ edges, with equality if and only if it is a forest, and the number of
edges is $\tau$.
\end{proof}

\begin{corollary}[Closure of two atoms]\label{cor:closure2}
Let $\rho,\rho'$ be atoms and let $D_1,\dots,D_l$ be the blocks of $\rho\vee\rho'$, in
order. Then $\mathrm{cl}(\{\rho,\rho'\})$ consists of the atoms obtained by
concatenating orderings of $D_1,\dots,D_l$, and it has $\prod_j|D_j|!$ elements. The
blocks are the classes of the equivalence relation on $[n]$ generated by the relation
``inverted'', that is, of its transitive closure.
\end{corollary}

\begin{proof}
The first statement is Lemma~\ref{lem:atoms} applied to $x=\rho\vee\rho'$. The second
is Theorem~\ref{thm:components}: the connected components of a graph are the classes of
the equivalence relation generated by its edges.
\end{proof}

\begin{example}\label{ex:delta}
Let $n=4$ and $\rho=(1\mid2\mid3\mid4)$. For $\rho'=(2\mid1\mid4\mid3)$ the graph has
the two edges $\{1,2\}$ and $\{3,4\}$, so $\delta=\tau=2$ and
$\rho\vee\rho'=(12\mid34)$, as in Example~\ref{ex:twopseudo}. For the reversal
$\rho'=(4\mid3\mid2\mid1)$ the graph is the complete graph $K_4$, so $\tau=6$ while
$\delta=3$ and $\rho\vee\rho'=\top$. For $\rho'=(3\mid4\mid1\mid2)$ the graph is the
$4$-cycle with edges $\{1,3\},\{1,4\},\{2,3\},\{2,4\}$, so $\tau=4$ and $\delta=3$:
the inequality is strict although the graph has no triangle.
\end{example}

The example with the reversal shows that $\delta$ saturates long before $\tau$: atoms
that are far apart in the Kendall sense have a join of large rank, but rank cannot
exceed $n-1$.

\begin{remark}[Forests and pattern avoidance]\label{rem:forest}
If $\rho$ is the identity, $\mathcal G(\rho,\rho')$ is a forest if and only if
$\rho'$ avoids the patterns $321$ and $3412$: a cycle in an inversion graph contains a
triangle (the pattern $321$) or an induced $4$-cycle (the pattern $3412$). We have also
checked this equivalence by computer for all permutations with $n\le8$. By Tenner
\cite[Theorem~4.3]{Tenner} these are exactly the permutations whose principal ideal in the
Bruhat order is a Boolean lattice, equivalently those with a reduced word without repeated
letters.
\end{remark}

\begin{remark}\label{rem:notcayley}
The distance $\delta$ should not be confused with the Cayley distance
$n-(\text{number of cycles of }\rho^{-1}\rho')$. For $\rho=(1\mid2\mid3\mid4)$ and
$\rho'=(3\mid4\mid1\mid2)$ the permutation $\rho^{-1}\rho'$ has two cycles, so the
Cayley distance is $2$, while $\delta=3$. Nor is it the Kendall distance, which is the
Hamming distance between the inversion sets, in the cube $\{0,1\}^{\binom n2}$: the
graph of the permutohedron is a partial cube for this embedding \cite{Ovchinnikov}.
\end{remark}

\section{The Möbius function and the characteristic polynomial}\label{sec:mobius}

In a geometric lattice the Möbius function $\mu(\bot,x)$ carries a lot of information:
its absolute value is not constant, the Whitney numbers of the first kind differ from
those of the second kind, and the characteristic polynomial is the main invariant.
Here the situation is the opposite. The Möbius function of $\mathcal L_n$ takes only the
values $\pm1$, with the sign dictated by the rank; as a consequence the Whitney
numbers of the two kinds coincide up to sign, and the characteristic polynomial is
determined by the numbers of elements of each rank. We show that these numbers
are the face numbers of the permutohedron, which brings in the Eulerian numbers,
and that the characteristic polynomial has a factor $t-1$ whose cofactor has
coefficients of strictly alternating sign, as in the geometric case. All
the proofs use the cut model.

\subsection{The Möbius function}

For a finite nonempty set $B$ let $\mathrm{OP}(B)$ be the set of ordered set
partitions of $B$, with the rank $\mathrm{rk}(z)=|B|-k(z)$ for $z\in\mathrm{OP}(B)$,
as in Definition~\ref{def:rank}.

\begin{lemma}[Elements below $x$]\label{lem:lowerset}
Let $x=(B_1\mid\cdots\mid B_k)\in\mathrm{OP}_n$. The elements $y\in\mathrm{OP}_n$ with
$y\le x$ are exactly the concatenations $y=y_1y_2\cdots y_k$ with
$y_j\in\mathrm{OP}(B_j)$, and the map $y\mapsto(y_1,\dots,y_k)$ is a bijection onto
$\prod_j\mathrm{OP}(B_j)$ with
\[
\mathrm{rk}(y)=\sum_{j=1}^k\mathrm{rk}(y_j).
\]
\end{lemma}

\begin{proof}
By definition $y\le x$ means that $x$ is obtained from $y$ by merging adjacent blocks,
that is, that each block $B_j$ of $x$ is the union of a collection of consecutive
blocks of $y$, and these collections follow one another in order. This says that $y$
is the concatenation of ordered set partitions $y_j$ of $B_j$, and the $y_j$ are
determined by $y$. If $y_j$ has $k_j$ blocks, then $y$ has $\sum_jk_j$ blocks and
$\mathrm{rk}(y)=n-\sum_jk_j=\sum_j(|B_j|-k_j)$.
\end{proof}

\begin{lemma}[An alternating count]\label{lem:alt}
For every finite nonempty set $B$,
\[
\sum_{z\in\mathrm{OP}(B)}(-1)^{\mathrm{rk}(z)}=1.
\]
\end{lemma}

\begin{proof}
Let $m=|B|$. The number of ordered set partitions of $B$ with $k$ blocks is
$k!\,S(m,k)$, and the exponential generating function of sequences of $k$ nonempty
sets is $(e^x-1)^k$, so that
\[
\sum_{m\ge0}\ \sum_k k!\,S(m,k)\,u^k\,\frac{x^m}{m!}=\frac1{1-u(e^x-1)}.
\]
For $u=-1$ the right-hand side is $e^{-x}$, and comparing the coefficients of $x^m/m!$
gives $\sum_k(-1)^kk!\,S(m,k)=(-1)^m$. Multiplying by $(-1)^m$ gives the claim,
since $\mathrm{rk}(z)=m-k$ for $z$ with $k$ blocks.
\end{proof}

\begin{theorem}[The Möbius function]\label{thm:mobius}
For all $x\le y$ in $\mathcal L_n$,
\[
\mu(x,y)=(-1)^{\mathrm{rk}(y)-\mathrm{rk}(x)}.
\]
In particular $\mu(\bot,x)=(-1)^{\mathrm{rk}(x)+1}$, and $|\mu(x,y)|=1$ for every
interval.
\end{theorem}

\begin{proof}
If $x\in\mathrm{OP}_n$, the interval $[x,y]$ is a Boolean lattice of rank
$\mathrm{rk}(y)-\mathrm{rk}(x)$ by Proposition~\ref{prop:boolean}, and the Möbius function
of a Boolean lattice is $(-1)^{\text{rank}}$. If $x=y=\bot$ the value is $1$. It
remains to treat $\mu(\bot,y)$ for $y\in\mathrm{OP}_n$, by induction on $\mathrm{rk}(y)$;
write $\mu(y)=\mu(\bot,y)$. Assume the formula for all $z<y$. By
Lemma~\ref{lem:lowerset} and Lemma~\ref{lem:alt},
\[
\sum_{z\in\mathrm{OP}_n,\ z\le y}(-1)^{\mathrm{rk}(z)}=\prod_{j=1}^k
\sum_{z_j\in\mathrm{OP}(B_j)}(-1)^{\mathrm{rk}(z_j)}=1,
\]
so that $\sum_{z\in\mathrm{OP}_n,\ z<y}(-1)^{\mathrm{rk}(z)}=1-(-1)^{\mathrm{rk}(y)}$. Hence
\[
\sum_{z<y}\mu(z)=\mu(\bot)+\sum_{z\in\mathrm{OP}_n,\ z<y}(-1)^{\mathrm{rk}(z)+1}
=1-\bigl(1-(-1)^{\mathrm{rk}(y)}\bigr)=(-1)^{\mathrm{rk}(y)},
\]
and the defining relation $\sum_{z\le y}\mu(z)=0$ gives $\mu(y)=(-1)^{\mathrm{rk}(y)+1}$.
\end{proof}

\begin{remark}[Eulerian lattices]\label{rem:eulerian}
A graded lattice with $\mu(x,y)=(-1)^{\mathrm{rk}(y)-\mathrm{rk}(x)}$ for all $x\le y$ is
called \emph{Eulerian}. By Remark~\ref{rem:geometry}, $\mathcal L_n$ is the face
lattice of the permutohedron, and face lattices of polytopes are Eulerian
\cite{StanleyEC1,Ziegler}; Theorem~\ref{thm:mobius} gives a direct proof in the cut model, which is
all we need. Thus $\mathcal L_n$ is pseudo-geometric and Eulerian, and its Möbius
function is as far as possible from that of a geometric lattice: for the
partition lattice $\Pi_n$, for instance, $|\mu(\hat0,\hat1)|=(n-1)!$.
\end{remark}

The value of $\mu$ also gives the inverse of Proposition~\ref{prop:atomsets}(3).

\begin{corollary}[Möbius inversion of the size]\label{cor:inversion}
For every $x\in\mathcal L_n$,
\[
\sum_{y\le x}(-1)^{\mathrm{rk}(x)-\mathrm{rk}(y)}\,s(y)=
\begin{cases}1&\text{if $x$ is an atom,}\\0&\text{otherwise.}\end{cases}
\]
\end{corollary}

\begin{proof}
By Proposition~\ref{prop:atomsets}(3), $s(x)=\sum_{y\le x}f(y)$ with $f$ the indicator
of the atoms, so $f(x)=\sum_{y\le x}\mu(y,x)s(y)$, and $\mu(y,x)$ is given by
Theorem~\ref{thm:mobius}.
\end{proof}

\subsection{Möbius inversion}\label{sec:inversion}

Since $\mathcal L_n$ is Eulerian, $\mu(x,y)=(-1)^{\operatorname{rk}y-\operatorname{rk}x}$ for all $x\le y$ (Theorem~\ref{thm:mobius}), and Möbius inversion has no coefficients other than signs. Write $k(x)$ for the number of blocks of $x\in\mathrm{OP}_n$; then $\operatorname{rk}y-\operatorname{rk}x=k(x)-k(y)$ for $x,y\ne\bot$.

\begin{proposition}\label{prop:inversion}
Let $F:\mathcal L_n\to\mathbb R$.
\begin{enumerate}[label=(\roman*)]
\item \emph{From below.} If $G(y)=\sum_{x\le y}F(x)$, then
\[
F(y)=\sum_{x\le y}(-1)^{\operatorname{rk}y-\operatorname{rk}x}\,G(x).
\]
The elements $x\le y$, $x\ne\bot$, are obtained by replacing each block of $y$ by an ordered set partition of that block, and $\bot$ contributes $(-1)^{\operatorname{rk}y+1}G(\bot)$ with $G(\bot)=F(\bot)$.
\item \emph{From above.} If $G(x)=\sum_{y\ge x}F(y)$, then
\[
F(x)=\sum_{S\subseteq C(x)}(-1)^{|S|}\,G(x\setminus S),
\]
where $x\setminus S$ is obtained from $x$ by deleting the cuts in $S$ (merging the corresponding adjacent blocks). This is the inclusion--exclusion formula on the Boolean interval $[x,\top]$.
\item \emph{An involution.} Let $T$ be the operator $(TF)(y)=(-1)^{\operatorname{rk}y}\sum_{x\le y}F(x)$. Then $T^2=\mathrm{id}$.
\item \emph{Multiplicative functions.} Let $F(x)=\prod_{B\in x}\varphi(B)$ for $x\ne\bot$ and $F(\bot)=c$, and put $\Phi(B)=\sum_{\pi}\prod_{C\in\pi}\varphi(C)$, the sum over all ordered set partitions $\pi$ of $B$. Then $G(y)=c+\prod_{B\in y}\Phi(B)$ for $y\ne\bot$, and the inversion of (i) gives back $\varphi(B)=\sum_{\pi}(-1)^{|\pi|-1}\prod_{C\in\pi}\Phi(C)$. If $\varphi(B)=\varphi_{|B|}$ and $\Phi(B)=\Phi_{|B|}$ depend only on sizes, with exponential generating functions $\varphi(t)=\sum_{m\ge1}\varphi_mt^m/m!$ and $\Phi(t)=\sum_{m\ge1}\Phi_mt^m/m!$, this is the pair of identities
\[
\Phi(t)=\frac{\varphi(t)}{1-\varphi(t)},\qquad \varphi(t)=\frac{\Phi(t)}{1+\Phi(t)}.
\]
\end{enumerate}
\end{proposition}

\begin{proof}
(i) and (ii) are Möbius inversion with $\mu(x,y)=(-1)^{\operatorname{rk}y-\operatorname{rk}x}$; in (ii) the interval $[x,\top]$ is Boolean on $C(x)$ and $\operatorname{rk}(x\setminus S)-\operatorname{rk}x=|S|$. (iii): if $Z$ is the zeta matrix and $D=\operatorname{diag}((-1)^{\operatorname{rk}})$, then $Z^{-1}=DZD$ by (i), so $(DZ)^2=DZDZ=Z^{-1}Z=I$. (iv): for $x\ne\bot$, $x\le y$ means that $x$ is a tuple of ordered set partitions $(\pi_B)_{B\in y}$, with $k(x)-k(y)=\sum_B(|\pi_B|-1)$ and $F(x)=\prod_B\prod_{C\in\pi_B}\varphi(C)$, so all sums factorize over the blocks of $y$; the constant $c$ cancels because $\sum_{x\le y}\mu(x,y)=0$. The generating functions follow from the fact that an ordered set partition is a sequence of blocks.
\end{proof}

\begin{example}
For $F=\mathbf 1_{\text{atoms}}$ we get $G(y)=\#\{\text{atoms below }y\}=\prod_j b_j!$, and (i) returns $1$ on atoms and $0$ elsewhere, as in Corollary~\ref{cor:inversion}. For $\varphi\equiv1$ one has $\Phi_m=$ the Fubini numbers and $\varphi(t)=e^t-1$, $\Phi(t)=1/(2-e^t)-1$.
\end{example}

\subsection{Whitney numbers}

For $-1\le k\le n-1$ let
\[
W_k=|\{x\in\mathcal L_n:\mathrm{rk}(x)=k\}|,\qquad
w_k=\sum_{\mathrm{rk}(x)=k}\mu(\bot,x),
\]
the Whitney numbers of the second and of the first kind. So $W_{-1}=1$ and $W_0=n!$.

\begin{theorem}[Whitney numbers]\label{thm:whitney}
For $0\le k\le n-1$,
\[
W_k=(n-k)!\,S(n,n-k),\qquad w_k=(-1)^{k+1}W_k,
\]
and $w_{-1}=W_{-1}=1$. Thus the Whitney numbers of the first kind are, up to sign, those
of the second kind: the numbers of faces of each dimension of the permutohedron.
\end{theorem}

\begin{proof}
An element of rank $k$ is an ordered set partition with $n-k$ blocks, and there are
$(n-k)!\,S(n,n-k)$ of them. The value of $w_k$ is Theorem~\ref{thm:mobius}.
\end{proof}

\begin{center}
\begin{tabular}{@{}c|cccccccc@{}}
$k$ & $-1$ & $0$ & $1$ & $2$ & $3$ & $4$ & $5$ & $|\mathcal L_n|$\\\hline
$n=3$ & $1$ & $6$ & $6$ & $1$ & & & & $14$\\
$n=4$ & $1$ & $24$ & $36$ & $14$ & $1$ & & & $76$\\
$n=5$ & $1$ & $120$ & $240$ & $150$ & $30$ & $1$ & & $542$\\
$n=6$ & $1$ & $720$ & $1800$ & $1560$ & $540$ & $62$ & $1$ & $4684$
\end{tabular}
\end{center}

\medskip
Let $A(n,i)$ be the Eulerian number, the number of permutations of $[n]$ with $i$
descents, and $A_n(x)=\sum_{i=0}^{n-1}A(n,i)x^i$ the Eulerian polynomial, so that
$A(n,0)=1$ and $A(n,i)=A(n,n-1-i)$.

\begin{proposition}[Whitney numbers and Eulerian numbers]\label{prop:eulerian}
\[
\sum_{k=0}^{n-1}W_ks^k=A_n(1+s)=\sum_{i=0}^{n-1}A(n,i)(1+s)^i.
\]
Consequently the polynomial $\sum_kW_ks^k$ has only real roots, all smaller than $-1$,
and the sequence $W_0,\dots,W_{n-1}$ is strictly log-concave, hence unimodal.
\end{proposition}

\begin{proof}
Since $W_k=(n-k)!\,S(n,n-k)$, the identity is Frobenius' formula
$A_n(x)=\sum_{m=1}^nm!\,S(n,m)(x-1)^{n-m}$ \cite{Frobenius,Petersen} with $x=1+s$ and $m=n-k$.
The Eulerian polynomial has $n-1$ simple negative roots \cite{Frobenius,Petersen}, so
$\sum_kW_ks^k=A_n(1+s)$ has $n-1$ simple real roots, each less than $-1$; Newton's
inequalities give the strict log-concavity of the coefficients, and a log-concave
sequence of positive numbers is unimodal.
\end{proof}

Since the permutohedron is a simple polytope, the identity says that the Eulerian numbers
form its $h$-vector; we do not use this.

\subsection{The characteristic polynomial}

\begin{definition}\label{def:chi}
The \emph{characteristic polynomial} of $\mathcal L_n$ is
\[
\chi_n(t)=\sum_{x\in\mathcal L_n}\mu(\bot,x)\,t^{\,\mathrm{rk}(\top)-\mathrm{rk}(x)}
=\sum_{x\in\mathcal L_n}\mu(\bot,x)\,t^{\,n-1-\mathrm{rk}(x)}.
\]
\end{definition}

Its degree is $n$, because $\mathrm{rk}(\bot)=-1$. By Theorem~\ref{thm:whitney},
\[
\chi_n(t)=t^n+\sum_{k=0}^{n-1}(-1)^{k+1}W_k\,t^{n-1-k}
=t^n-n!\,t^{n-1}+\tfrac{(n-1)\,n!}2\,t^{n-2}-\cdots+(-1)^n.
\]
The coefficient of $t^{n-2}$ is $W_1=(n-1)\,n!/2$, the number of pseudo-atoms.

\begin{theorem}[Structure of $\chi_n$]\label{thm:chi}
Let $n\ge3$.
\begin{enumerate}
\item $\chi_n(t)=t^n-\sum_{i=0}^{n-1}A(n,i)\,(t-1)^i\,t^{\,n-1-i}$.
\item $\chi_n(t)=(t-1)\,\bar\chi_n(t)$ with
  \[
  \bar\chi_n(t)=\sum_{j=0}^{n-1}(-1)^j\beta_j\,t^{\,n-1-j},
  \]
  where $\beta_0=1$ and $\beta_j=\sum_{i=j}^{n-1}A(n,i)\binom{i-1}{j-1}$ for $j\ge1$;
  all $\beta_j$ are positive integers. In particular
  $\beta_1=n!-1$, $\beta_{n-2}=2^n-3$ and $\beta_{n-1}=1$.
\item All real roots of $\chi_n$ are positive, and $\chi_n(0)=(-1)^n$,
  $\bar\chi_n(1)=\chi_n'(1)=n+2-2^n$.
\end{enumerate}
\end{theorem}

\begin{proof}
Put $R_n(w)=\sum_{x\in\mathcal L_n}w^{\mathrm{rk}(x)+1}=1+\sum_{k=0}^{n-1}W_kw^{k+1}$, the
rank generating polynomial of $\mathcal L_n$. By Theorem~\ref{thm:mobius},
$\mu(\bot,x)\,t^{\,n-1-\mathrm{rk}(x)}=(-1)^{\mathrm{rk}(x)+1}t^{\,n-(\mathrm{rk}(x)+1)}$, so
\[
\chi_n(t)=t^n\,R_n(-1/t).
\]
By Proposition~\ref{prop:eulerian} and $A(n,0)=1$,
\[
R_n(w)=1+w\sum_{i=0}^{n-1}A(n,i)(1+w)^i=(1+w)\Bigl(1+w\sum_{i=1}^{n-1}A(n,i)(1+w)^{i-1}\Bigr).
\]
(1) Substitute $w=-1/t$, so that $1+w=(t-1)/t$:
$\chi_n(t)=t^n-\sum_iA(n,i)(t-1)^it^{n-1-i}$.

\noindent (2) Expanding $w(1+w)^{i-1}=\sum_{j=1}^{i}\binom{i-1}{j-1}w^j$ gives
$R_n(w)=(1+w)\sum_{j=0}^{n-1}\beta_jw^j$ with the $\beta_j$ of the statement. Then
\[
\chi_n(t)=t^n\Bigl(1-\frac1t\Bigr)\sum_j\beta_j\Bigl(-\frac1t\Bigr)^j
=(t-1)\sum_j(-1)^j\beta_jt^{n-1-j}.
\]
The $\beta_j$ are positive since $A(n,i)>0$ for $0\le i\le n-1$. Moreover
$\beta_1=\sum_{i\ge1}A(n,i)=n!-1$ and $\beta_{n-1}=A(n,n-1)=1$. Finally
$\beta_{n-2}=A(n,n-2)+(n-2)\,A(n,n-1)=A(n,1)+n-2=2^n-3$, because $A(n,n-2)=A(n,1)$ and
$A(n,1)=2^n-n-1$.
(3) For $t\le0$ every term $(-1)^j\beta_jt^{n-1-j}$ of $\bar\chi_n$ has the sign of
$(-1)^{n-1}$, so $\bar\chi_n$ has no root in $(-\infty,0]$, and neither has $\chi_n$.
The value $\chi_n(0)=\mu(\bot,\top)=(-1)^n$ is Theorem~\ref{thm:mobius}. Finally
$\chi_n'(1)=\bar\chi_n(1)$ and, differentiating (1) at $t=1$, only the terms with $i=0$ and
$i=1$ contribute: $\chi_n'(1)=n-(n-1)-A(n,1)=n+2-2^n$.
\end{proof}

Part (2) is the analogue of the sign alternation of the coefficients of
$\chi_L(t)/(t-1)$ for a geometric lattice $L$ \cite{Rota64}. In the geometric case the
absolute values of these coefficients are the Betti numbers of the broken circuit
complex. Here we do not have a complex, but the $\beta_j$ do count something explicit.
Write permutations of $[n]$ as words $\sigma=\sigma_1\cdots\sigma_n$, with
$\mathrm{Des}(\sigma)=\{i:\sigma_i>\sigma_{i+1}\}$.

\begin{proposition}[Combinatorial meaning of $\beta_j$]\label{prop:betaj}
\begin{enumerate}
\item For $1\le j\le n-1$, $\beta_j$ is the number of pairs $(\sigma,M)$ with
  $\sigma\in S_n$, $M\subseteq\mathrm{Des}(\sigma)$, $|M|=j$ and $\min\mathrm{Des}(\sigma)\in M$.
\item There is a perfect matching of the Hasse diagram of $\mathcal L_n$ in which, for
  every $0\le j\le n-1$, exactly $\beta_j$ matched pairs $x\lessdot y$ have
  $\mathrm{rk}(y)=j$. The matching is acyclic: orienting its edges downwards and all
  other cover relations upwards gives a digraph without directed cycles.
\end{enumerate}
\end{proposition}

\begin{proof}
(1) A permutation with $i\ge1$ descents has $\binom{i-1}{j-1}$ subsets of $j$ descents
containing the first one, so the number in the statement is
$\sum_iA(n,i)\binom{i-1}{j-1}=\beta_j$ by Theorem~\ref{thm:chi}(2).

(2) For $x=(B_1\mid\dots\mid B_k)\in\mathrm{OP}_n$ let $\sigma(x)$ be the word obtained by
listing the elements of each block in decreasing order, block after block, and let
$S(x)$ be the set of positions $i$ such that $\sigma(x)_i,\sigma(x)_{i+1}$ lie in the same
block. Then $S(x)\subseteq\mathrm{Des}(\sigma(x))$, $|S(x)|=n-k=\mathrm{rk}(x)$, and
$x\mapsto(\sigma(x),S(x))$ is a bijection onto the pairs $(\sigma,S)$ with
$S\subseteq\mathrm{Des}(\sigma)$ (cut $\sigma$ at the positions not in $S$).
Let $x_0=(1\mid2\mid\cdots\mid n)$, the only $x$ with $\sigma(x)$ having no descent,
and match $\bot$ with $x_0$. For $x\ne x_0$ let $d=\min\mathrm{Des}(\sigma(x))$ and let
$m(x)$ correspond to $(\sigma(x),S(x)\,\triangle\,\{d\})$. If $d\notin S(x)$, then
$m(x)$ is obtained from $x$ by merging the two adjacent blocks that contain the
positions $d$ and $d+1$ (the merged block is again decreasing in $\sigma(x)$), so
$x\lessdot m(x)$; if $d\in S(x)$ the roles are exchanged. Thus $m$ is a fixed-point free
involution on $\mathrm{OP}_n\setminus\{x_0\}$ whose pairs are cover relations, and
$\{\bot,x_0\}$ is one more pair. The pairs with $\mathrm{rk}(y)=j\ge1$ correspond to the
pairs $(\sigma,S')$ with $S'\subseteq\mathrm{Des}(\sigma)\setminus\{\min\mathrm{Des}(\sigma)\}$ and
$|S'|=j-1$, whose number is $\sum_\sigma\binom{\mathrm{des}(\sigma)-1}{j-1}=\beta_j$ as in (1);
for $j=0$ there is the pair $\{\bot,x_0\}$, and $\beta_0=1$.

Acyclicity. For $x\in\mathrm{OP}_n$ put $\Phi(x)=\sum_ii\,\sigma(x)_i$, and
$\Phi(\bot)=\sum_ii^2$. By the rearrangement inequality $\Phi(x)$ is the
\emph{unique} minimum of $\sum_ii\,w_i$ over the words $w$ obtained by ordering each block
of $x$ arbitrarily, and $\sum_ii^2$ is the unique maximum over all words. If $x<y$ in
$\mathrm{OP}_n$, the words of $y$ include those of $x$, so $\Phi(y)\le\Phi(x)$, with equality
only if $\sigma(y)=\sigma(x)$; and $\Phi(y)<\Phi(\bot)$ for every $y\ne x_0$, while
$\Phi(x_0)=\Phi(\bot)$.
Consider a directed cycle. Each upward step is a cover, so it does not increase $\Phi$, and
each downward step is a matched pair, which has the same $\Phi$; hence $\Phi$ is constant on the
cycle. The cycle does not contain $\bot$: it could only leave $\bot$ by an upward step to an atom
$y\ne x_0$ (the pair $\{\bot,x_0\}$ is matched), and then $\Phi(y)<\Phi(\bot)$. So
all elements of the cycle have the same $\sigma$, that is, they are pairs $(\sigma,S)$ with
$S\subseteq\mathrm{Des}(\sigma)$, and $d=\min\mathrm{Des}(\sigma)$ is fixed. A downward step removes
$d$ from $S$ and an upward step adds an element of $\mathrm{Des}(\sigma)$ different from $d$.
If the cycle has no downward step, $|S|$ increases at every step, which is impossible;
if it has one, then $d\in S$ just before that step and $d\notin S$ just after it, and no
step ever adds $d$; going once around the cycle we would need $d\in S$ again before the same step,
a contradiction.
\end{proof}

For $n=3$ the numbers $\beta_j=1,5,1$ are: the pair $\{\bot,x_0\}$; the five pairs
(vertex, edge) with vertex $\sigma\ne123$; and the pair (edge, hexagon) for $\sigma=321$.
We checked (1) by computer for $n\le8$, and (2), including acyclicity, for $n\le5$.

\begin{remark}[Toric reading]\label{rem:toric}
By Proposition~\ref{prop:eulerian}, $A_n(q)=\sum_F(q-1)^{\dim F}$, the sum over the
non-empty faces $F$ of the permutohedron; this is the number of points over $\mathbb F_q$
of the toric variety of the permutohedron, whose Betti numbers are the Eulerian numbers.
Likewise $R_n(q-1)=1+(q-1)A_n(q)$ is the number of points of the normal affine toric
variety of the cone over the permutohedron, one orbit for each element of
$\mathcal L_n$, and Theorem~\ref{thm:chi}(2) says that this number equals
$q\sum_j\beta_j(q-1)^j$. This explains the factor $t-1$ of $\chi_n$, but we know no
space whose Betti numbers are the $\beta_j$ themselves.
\end{remark}

\begin{example}
For $n=3$, $\chi_3(t)=t^3-6t^2+6t-1=(t-1)(t^2-5t+1)$, with roots
$1$ and $(5\pm\sqrt{21})/2$.

\noindent For $n=4$,
$\chi_4(t)=t^4-24t^3+36t^2-14t+1=(t-1)(t^3-23t^2+13t-1)$.

\noindent For $n=5$,
$\chi_5(t)=t^5-120t^4+240t^3-150t^2+30t-1=(t-1)(t^4-119t^3+121t^2-29t+1)$.
\end{example}

\begin{remark}[Beta invariant]
For a geometric lattice $L$, $|\chi_L'(1)|$ is Crapo's beta invariant \cite{Crapo}. By
part (3), the analogous quantity for $\mathcal L_n$ is
$|\chi_n'(1)|=2^n-n-2=A(n,1)-1$, one less than the number of permutations of $[n]$ with
exactly one descent.
\end{remark}

\subsection{Real roots}

By Theorem~\ref{thm:chi}(3) the real roots of $\chi_n$ are positive, and $t=1$ is a
root. A computation suggests that there are no others.

\begin{conjecture}[Real-rootedness]\label{conj:real}
For every $n\ge3$ the polynomial $\chi_n(t)$ has $n$ distinct positive real roots:
the root $1$, exactly $n-2$ roots in $(0,1)$ and one root greater than $1$.
Equivalently, the rank generating polynomial $R_n(w)=\sum_{x\in\mathcal L_n}w^{\mathrm{rk}(x)+1}$
has only real roots, all negative.
\end{conjecture}

We checked exactly, with Sturm sequences over the rationals, that $\chi_n$ has $n$
distinct real roots for $3\le n\le30$ and for $n=31,34,37,40,43$, and that for
$3\le n\le25$ the polynomial $\bar\chi_n$ has $n-2$ roots in $(0,1)$ and one in
$(1,n!+1)$. The sum of the roots is $n!$ and their product is $1$. With $u=1-1/t$,
the equation $\chi_n(t)=0$ becomes $(1-u)A_n(u)=1$: the conjecture says that the
graph of $(1-u)A_n(u)$ crosses the level $1$ in $n$ points, although the roots
of $A_n$ themselves are real by Proposition~\ref{prop:eulerian}. Real-rootedness
is not a property of characteristic polynomials of geometric lattices in general,
where one has instead the log-concavity of the coefficients \cite{AHK}.
The conjecture reduces to an inequality about the roots of the Eulerian polynomial
$A_{n-1}$, which are simple, real and negative \cite{Frobenius,Petersen}.

\begin{proposition}[Reduction of the conjecture]\label{prop:reduction}
Let $n\ge3$ and let $r_0>r_1>\dots>r_{n-3}$ be the roots of $A_{n-1}$. If
\[
|r_k|\,(1+|r_k|)^2\,|A_{n-1}'(r_k)|>1\qquad\text{for every odd }k,\tag{$\ast$}
\]
then Conjecture~\ref{conj:real} holds for $n$.
\end{proposition}

\begin{proof}
Put $P_n(u)=1+(u-1)A_n(u)$. By the proof of Theorem~\ref{thm:chi}, $P_n(u)=R_n(u-1)$ and
$\chi_n(t)=t^nP_n(1-1/t)$. The map $t\mapsto u=1-1/t$ is an increasing bijection of
$(0,\infty)$ onto $(-\infty,1)$ with $1\mapsto0$, $(0,1)\mapsto(-\infty,0)$ and
$(1,\infty)\mapsto(0,1)$. So it suffices to show that $P_n$ has $n-2$ distinct roots
in $(-\infty,0)$, the root $0$, and one root in $(0,1)$; then $\chi_n$ has the
root $1$, $n-2$ roots in $(0,1)$ and one root greater than $1$.

We have $P_n(0)=1-A(n,0)=0$, $P_n(1)=1$, and
$P_n'(0)=1-A(n,1)=-(2^n-n-2)<0$, so $P_n<0$ just to the right of $0$ and $P_n$ has a root
in $(0,1)$. Next, the recurrence $A_n(u)=(1+(n-1)u)A_{n-1}(u)+u(1-u)A_{n-1}'(u)$
gives, at a root $r$ of $A_{n-1}$, $A_n(r)=r(1-r)A_{n-1}'(r)$ and hence
\[
P_n(r)=1-r(1-r)^2A_{n-1}'(r).
\]
Since $A_{n-1}$ is monic with simple roots, $A_{n-1}'(r_k)$ has the sign $(-1)^k$. For
$k$ even, $P_n(r_k)=1+|r_k|(1+|r_k|)^2|A_{n-1}'(r_k)|>0$; for $k$ odd,
$P_n(r_k)=1-|r_k|(1+|r_k|)^2|A_{n-1}'(r_k)|<0$ by $(\ast)$. Therefore $P_n$ has a root in
each interval $(r_{k+1},r_k)$, which gives $n-3$ negative roots. Finally $P_n(u)$ has the
sign of $(-1)^n$ for $u\to-\infty$ and $P_n(r_{n-3})$ has the sign of $(-1)^{n-3}$, so
there is one more root in $(-\infty,r_{n-3})$. The roots found are $n-2$ negative ones,
$0$, and one in $(0,1)$: $n$ distinct roots in all, hence all of them.
\end{proof}

For $n=3$ the inequality is empty, so $\chi_3$ is real-rooted, as seen directly. We
computed the roots of $A_{n-1}$ to high precision for $4\le n\le27$ and found that the
left-hand side of $(\ast)$ is at least $11$ for every odd $k$, with minimum $11.07$ at
$n=8$ and increasing for larger $n$; the signs of $P_n(r_k)$ were also checked
directly. So the margin is large, but $(\ast)$ is a quantitative statement on the roots of
$A_{n-1}$ which we have not proved.

\section{Modular elements and supersolvability}\label{sec:modular}

In a geometric lattice, a maximal chain of modular elements, that is, supersolvability,
is the source of the factorization of the characteristic polynomial into linear factors
with integer roots, and of an easy EL-labelling \cite{Stanley72}. Here nothing of the kind
is available: we show that $\mathcal L_n$ has no left-modular element other than
$\bot$ and $\top$, and in particular no left-modular atom.

\subsection{Left-modular elements}

\begin{definition}\label{def:leftmodular}
An element $x$ of a lattice $L$ is \emph{left-modular} if
\[
y\vee(x\wedge z)=(y\vee x)\wedge z\qquad\text{for all }y\le z\text{ in }L.
\]
\end{definition}

For $y\le z$ the inequality $y\vee(x\wedge z)\le(y\vee x)\wedge z$ always holds, so
only the opposite inequality is at stake. The elements $\bot$ and $\top$ are
left-modular in every lattice: for $x=\top$ both sides equal $z$, and for $x=\bot$
both sides equal $y$.

A finite lattice is \emph{supersolvable} if it has a maximal chain $M$ such that, for
every chain $C$, the sublattice generated by $C\cup M$ is distributive \cite{Stanley72}.
Every element $x$ of such a chain $M$ is left-modular: for $y\le z$, the sublattice
generated by $x,y,z$ is distributive, so
$(y\vee x)\wedge z=(y\wedge z)\vee(x\wedge z)=y\vee(x\wedge z)$. Hence a
supersolvable lattice has a maximal chain of left-modular elements, which is also
the usual alternative definition, and it suffices to show that this fails for
$\mathcal L_n$.

\begin{theorem}[Left-modular elements]\label{thm:leftmodular}
Let $n\ge3$. The left-modular elements of $\mathcal L_n$ are exactly $\bot$ and $\top$.
\end{theorem}

\begin{proof}
It remains to show that no $x\in\mathrm{OP}_n$ with $x\ne\top$ is left-modular. Then
$C(x)$ is a nonempty chain of proper subsets of $[n]$. We construct $y\le z$ with
$y\vee(x\wedge z)\ne(y\vee x)\wedge z$.

Fix $A\in C(x)$, with $|A|=j$, and let $B\ne A$ be another subset with $|B|=j$. Two
distinct subsets of the same size are not comparable, so $B$ is not comparable with
$A$, and therefore $B\notin C(x)$, since $C(x)$ is a chain. Next we choose a cut $b$
comparable with $B$ and different from it, with $b\notin C(x)$. If $|B|\ge2$, the
$|B|$ subsets of $B$ of size $|B|-1$ are pairwise incomparable, so at most one of them
belongs to $C(x)$; take $b$ among the others. If $|B|=1$, then $n\ge3$, and the
$n-1\ge2$ sets $B\cup\{e\}$ with $e\notin B$ are pairwise incomparable; take $b$ among
those not in $C(x)$.

Since $B$ and $b$ are comparable, there is a flag $F$, that is, a maximal chain of
proper nonempty subsets of $[n]$, with $B,b\in F$. Let $y$ be the atom with
$C(y)=F$, and let $z$ be the element with $C(z)=F\setminus\{b\}$; it covers $y$ by
Proposition~\ref{prop:boolean}, so $y\le z$.

The set $C(x)\cup C(z)$ contains $A$ and $B$, which are not comparable, so it is not a
chain and $x\wedge z=\bot$ by Theorem~\ref{thm:lattice}(2). Hence
$y\vee(x\wedge z)=y$. On the other hand $C(y\vee x)=F\cap C(x)$ by
Theorem~\ref{thm:lattice}(1), and $b\notin C(x)$, so $F\cap C(x)\subseteq F\setminus\{b\}
=C(z)$. Thus $C(y\vee x)\cup C(z)=C(z)$ is a chain, and
$(y\vee x)\wedge z=z$ by Theorem~\ref{thm:lattice}(2). Since $z\ne y$, the element
$x$ is not left-modular.
\end{proof}

\begin{example}
Let $n=3$, $x=(1\mid2\mid3)$ with $C(x)=\{1,12\}$. Take $A=1$, $B=3$ and $b=32$; the
flag $F=\{3,32\}$ gives $y=(3\mid2\mid1)$ and $z=(3\mid12)$. Then $x\wedge z=\bot$ since
$\{1,12,3\}$ is not a chain, while $y\vee x=\top$ since $F\cap C(x)=\emptyset$, and
$\top\wedge z=z\ne y$.
\end{example}

\begin{corollary}[No supersolvability]\label{cor:notss}
For $n\ge3$ no atom of $\mathcal L_n$ is left-modular. Hence $\mathcal L_n$ has no
maximal chain of left-modular elements and is not supersolvable.
\end{corollary}

\begin{proof}
Every maximal chain contains an atom, and the only left-modular elements are $\bot$ and
$\top$ by Theorem~\ref{thm:leftmodular}.
\end{proof}

\subsection{No factorization of the characteristic polynomial}

For a supersolvable semimodular lattice the characteristic polynomial is a product of
linear factors $t-a_i$ with nonnegative integers $a_i$ \cite{Stanley72}. The polynomial
$\chi_n$ of Section~\ref{sec:mobius} has no such factorization, independently of
Theorem~\ref{thm:leftmodular}.

\begin{proposition}[Roots of $\chi_n$]\label{prop:noroots}
For $n\ge3$ the only rational root of $\chi_n$ is the simple root $1$. In particular
$\chi_n$ is not a product of linear factors with integer coefficients.
\end{proposition}

\begin{proof}
The polynomial $\chi_n$ is monic with integer coefficients and constant term
$(-1)^n$, so its rational roots are $\pm1$. All its real roots are positive by
Theorem~\ref{thm:chi}(3), and $1$ is a simple root because
$\bar\chi_n(1)=n+2-2^n\ne0$ for $n\ge3$. Since $\deg\bar\chi_n=n-1\ge2$, the
cofactor $\bar\chi_n$ has no rational root and is not a product of linear factors over
$\mathbb Z$.
\end{proof}

For $n=3$ the roots are $1$ and $(5\pm\sqrt{21})/2$, which are irrational.

\begin{remark}
The Eulerian property and the alternation of signs of Section~\ref{sec:mobius}, which
make $\mathcal L_n$ resemble a geometric lattice, therefore do not come from modular
chains. Whether $\mathcal L_n$ is EL-shellable is a separate question, taken up in
Section~\ref{sec:el}.
\end{remark}

\section{Shellability}\label{sec:shell}

Shellability links the order structure of a graded poset with the topology of its order
complex. A geometric lattice is EL-shellable \cite{Bjorner80}. Here the lattice is the face
lattice of a simple polytope, and what is inherited from the polytope is the weaker notion
of CL-shellability \cite{BjornerWachs}. We recall the notions we need, show that every
linear extension of the weak order on permutations is a recursive atom ordering of
$\mathcal L_n$ (this contains the ordering used by Sagan and Sundaram
\cite{SaganSundaram}), record what follows for the flag $h$-vector, and then discuss what is
missing for EL-shellability.

\subsection{Definitions}

Let $P$ be a bounded graded poset with minimum $\hat0$ and maximum $\hat1$. A
\emph{recursive atom ordering} of $P$ is defined by induction on the length of $P$
\cite[Definition~3.1]{BjornerWachs}. If $P$ has length $1$, every ordering of its (at most
one) atom is one. If the length is at least $2$, a total order $a_1,\dots,a_t$ of the atoms
of $P$ is a recursive atom ordering if

\begin{enumerate}[label=(R\arabic*)]
\item for every $j$ the interval $[a_j,\hat1]$ has a recursive atom ordering in which the
  atoms of $[a_j,\hat1]$ that cover some $a_i$ with $i<j$ come first, and
\item for $i<j$, if $a_i$ and $a_j$ are both below an element $y$, then there are $k<j$ and
  an atom $z$ of $[a_j,\hat1]$ such that $a_k<z\le y$.
\end{enumerate}

A \emph{CL-labelling} of $P$ assigns to each cover relation of each rooted interval (an
interval $[x,y]$ together with a maximal chain of $[\hat0,x]$) a label in a poset, in such a
way that every rooted interval has a unique maximal chain with increasing labels, and this
chain is lexicographically smaller than every other maximal chain of the interval; an
EL-labelling is a CL-labelling in which the labels of a cover relation do not depend on
the chosen root. The poset $P$ is \emph{CL-shellable} (respectively \emph{EL-shellable}) if
it has such a labelling.

\begin{theorem}[Bj\"orner--Wachs]\label{thm:BW}
A bounded graded poset is CL-shellable if and only if it has a recursive atom ordering
\cite[Theorem~3.2]{BjornerWachs}. If $P$ is CL-shellable then every interval of $P$ is
CL-shellable, the order complex of the proper part $\bar P=P\setminus\{\hat0,\hat1\}$ is
shellable, hence Cohen--Macaulay, and it is homotopy equivalent to a wedge of
$|\mu_P(\hat0,\hat1)|$ spheres of dimension $\mathrm{length}(P)-2$. Moreover, for a
CL-labelling and a set $S$ of ranks, the flag $h$-number $\beta_P(S)$ is the number of
maximal chains whose labels descend exactly at the positions in $S$ \cite[Theorem~10.6]{BjornerWachsII}.
\end{theorem}

Here the flag $h$-numbers are defined from the flag $f$-numbers $\alpha_P(T)$, the numbers
of chains of $\bar P$ whose elements have ranks in $T$, by
$\beta_P(S)=\sum_{T\subseteq S}(-1)^{|S\setminus T|}\alpha_P(T)$. We use the rank shifted by
one, so that the elements of rank $r\ge0$ of $\mathcal L_n$ have position $r+1\in[n-1]$,
and $S$ ranges over the subsets of $[n-1]$.

\subsection{Recursive atom orderings from the weak order}

We write an atom as a word $\rho=\rho_1\rho_2\cdots\rho_n$, the permutation
$(\rho_1\mid\cdots\mid\rho_n)$. The (right) weak order on permutations is
defined by $u\preceq w$ if every pair of values $a<b$ that is inverted in $u$ (that is, $b$
stands before $a$) is also inverted in $w$. If $w_p>w_{p+1}$ at two adjacent positions,
exchanging these two letters produces a word $u$ with $u\prec w$, whose inversions are those
of $w$ except for the pair $\{w_{p+1},w_p\}$. In particular the lexicographic order of words
is a linear extension of $\preceq$, since such an exchange makes the word lexicographically
smaller.

\begin{theorem}[Recursive atom orderings of $\mathcal L_n$]\label{thm:RAO}
Let $n\ge2$ and let $\pi$ be any total order of the permutations of $[n]$ that extends the
weak order, that is, $u\preceq w$ implies that $u$ precedes $w$ in $\pi$. Then $\pi$ is a
recursive atom ordering of $\mathcal L_n$.
\end{theorem}

\begin{proof}
First, in a Boolean lattice every total order of the atoms is a recursive atom ordering.
Indeed, for (R2) let $a_i,a_j$ be atoms below $y$ with $i<j$, and take $k=i$ and
$z=a_i\vee a_j$, which is an atom of $[a_j,\hat1]$ and satisfies $a_i<z\le y$. For (R1), the
interval above an atom is again Boolean, so by induction on the rank it has a recursive atom
ordering for every total order of its atoms, in particular for one in which prescribed atoms
come first. Now (R1) holds for $\pi$, because the interval $[\rho,\top]$ above an atom $\rho$
of $\mathcal L_n$ is Boolean by Proposition~\ref{prop:boolean}.

It remains to prove (R2) for $\pi$. Let $\rho,\sigma$ be atoms with $\sigma$ preceding
$\rho$ in $\pi$, both below some $y=(B_1\mid\cdots\mid B_k)\in\mathrm{OP}_n$. The atoms below
$y$ are the words whose first $|B_1|$ letters form $B_1$, whose next $|B_2|$ letters form
$B_2$, and so on, and this set is closed under exchanging two adjacent letters that lie in
the same block. Let $\rho^+$ be the word obtained by writing the letters of each block in
increasing order. Repeated exchanges of adjacent descents inside the blocks lead from any
atom $w\le y$ to $\rho^+$, so $\rho^+\preceq w$, and $\rho^+$ precedes every other atom
below $y$ in $\pi$. Since $\sigma$ precedes $\rho$, the atom $\rho$ is not the $\pi$-first
atom below $y$, so $\rho\ne\rho^+$, and $\rho$ has two adjacent letters $\rho_p>\rho_{p+1}$
in the same block of $y$. Exchanging them gives an atom $\rho'\prec\rho$ below $y$, which
precedes $\rho$ in $\pi$. The element $z$ of $\mathrm{OP}_n$ obtained from $\rho$ by merging
the singletons $\{\rho_p\}$ and $\{\rho_{p+1}\}$ into one block covers both $\rho$ and
$\rho'$, so it is an atom of $[\rho,\top]$ with $\rho'<z$, and $z\le y$ because the merged
block lies inside a block of $y$. This is (R2) with $a_k=\rho'$.
\end{proof}

\begin{corollary}[CL-shellability]\label{cor:CL}
$\mathcal L_n$ is CL-shellable for every $n\ge2$. The order complex of its proper part
is shellable and homotopy equivalent to a sphere of dimension $n-2$. The same holds for every
interval of $\mathcal L_n$.
\end{corollary}

\begin{proof}
Theorem~\ref{thm:RAO} and Theorem~\ref{thm:BW}. By Theorem~\ref{thm:mobius},
$|\mu(\bot,\top)|=1$ and the length of $\mathcal L_n$ is $n$, so the proper part
has the homotopy type of a sphere of dimension $n-2$.
\end{proof}

\begin{remark}\label{rem:inherited}
The lexicographic order of the permutations is a linear extension of the weak order, so
Theorem~\ref{thm:RAO} contains the recursive atom ordering found by Sagan and
Sundaram \cite[Theorem~2.7]{SaganSundaram}, and any other linear extension works as well.
Corollary~\ref{cor:CL} is also a special case of the fact that the face lattice of every
polytope is CL-shellable \cite[Theorem~4.5]{BjornerWachs}, which is where it comes from:
line shellings of polytopes \cite[Chapter~8]{Ziegler}, applied to the polar of the
permutohedron, which is the Coxeter complex of type $A_{n-1}$. We claim no novelty for the
shellability itself; what Theorem~\ref{thm:RAO} adds is the explicit description of all the
atom orderings that arise from the weak order.
\end{remark}

\subsection{The $h$-vector and the flag $h$-vector}

The order complex of the proper part of $\mathcal L_n$ is the barycentric subdivision of
the boundary complex of the polar of the permutohedron, that is, of the Coxeter complex of
type $A_{n-1}$, a simplicial sphere of dimension $n-2$ whose $h$-vector is the sequence of
Eulerian numbers $A(n,0),\dots,A(n,n-1)$. Brenti and Welker computed the $h$-vector of the
barycentric subdivision of a simplicial complex $\Delta$ of dimension $d-1$ from the
$h$-vector of $\Delta$: it is $h_j(\mathrm{sd}\,\Delta)=\sum_{r=0}^dA(d+1,j,r+1)\,h_r(\Delta)$,
where $A(d,j,i)$ is the number of permutations of $[d]$ with $j$ descents and first letter
$i$ \cite[Theorem~2.2]{BrentiWelker}. The refinement by descent sets is a theorem of Stanley
\cite[Theorem~8.3]{Stanley82}, in the form given by Sagan and Sundaram
\cite[Theorem~4.2 and Proposition~4.10]{SaganSundaram}; applied to the Boolean lattice it gives the
flag $h$-numbers of $\mathcal L_n$, as follows.

\begin{theorem}[Flag $h$-numbers of $\mathcal L_n$]\label{prop:flagh}
For $n\ge3$ and every $S\subseteq[n-1]$,
\[
\beta_{\mathcal L_n}(S)=\sum_{\substack{w\in\mathfrak S_n\\ \mathrm{Des}(w)=S}}A(n,w_1-1),
\]
where $A(n,i)$ is the Eulerian number and $w_1$ is the first letter of $w$. Summing over
$|S|=j$, with $A(n,j,i)$ denoting the number of
permutations of $[n]$ with $j$ descents and first letter $i$, gives $h_j=\sum_{r=0}^{n-1}A(n,j,r+1)\,A(n,r)$; this is the formula of Brenti and
Welker for the Coxeter complex.
\end{theorem}

\begin{proof}
Let $B_n$ be the Boolean lattice of the subsets of $[n]$, of rank $n$, and let $B_n^{\mathrm{sd}}$ be
the poset of chains of the proper part of $B_n$, ordered by inclusion, with a top element adjoined
(the empty chain is the minimum). By the cut model, $x\mapsto C(x)$ is an order-reversing bijection
from $\mathrm{OP}_n$ onto the chains of proper nonempty subsets (the empty chain included), so $\mathcal L_n$ is the dual
of $B_n^{\mathrm{sd}}$; an element $x$ with $|C(x)|=c$ has rank $c$ in $B_n^{\mathrm{sd}}$ and rank
$n-c$ in $\mathcal L_n$, when ranks are shifted by one as in this section. The order complex of the
proper part is the same for both, and rank selection by $T$ in $B_n^{\mathrm{sd}}$ is rank selection by
$n-T=\{n-t:t\in T\}$ in $\mathcal L_n$, so $\beta_{\mathcal L_n}(S)=\beta_{B_n^{\mathrm{sd}}}(n-S)$.

Stanley's theorem, in the formulation of \cite[Theorem~4.2]{SaganSundaram}, says that for a
Cohen--Macaulay bounded graded poset $P$ of rank $n$,
$\beta_{P^{\mathrm{sd}}}(T)=\sum_{i=0}^{n-1}c_{i,n}(T)\,\eta_i(P)$ in dimension, where $c_{i,n}(T)$
is the number of $\sigma\in\mathfrak S_n$ with $\sigma(n)=n-i$ and $\mathrm{Des}(\sigma)=T$, and
$\eta_i(P)=\sum_{|U|=i}\beta_P(U)$. For $P=B_n$ one has $\beta_{B_n}(U)=\#\{\pi\in\mathfrak S_n:
\mathrm{Des}(\pi)=U\}$ (Solomon, Stanley; \cite[Theorem~4.1]{SaganSundaram}), so $\eta_i(B_n)=A(n,i)$. Hence
\[
\beta_{\mathcal L_n}(S)=\sum_{\substack{\sigma\in\mathfrak S_n\\ \mathrm{Des}(\sigma)=n-S}}A(n,\,n-\sigma(n)).
\]
Finally put $w(j)=n+1-\sigma(n+1-j)$. Then $\sigma\mapsto w$ is a bijection of $\mathfrak S_n$,
$\mathrm{Des}(w)=\{n-t:t\in\mathrm{Des}(\sigma)\}$, which is $S$ when $\mathrm{Des}(\sigma)=n-S$, and
$w_1=n+1-\sigma(n)$, so $A(n,n-\sigma(n))=A(n,w_1-1)$. The second formula follows by grouping the
permutations $w$ by the number of their descents and by their first letter, since $A(n,r)$ depends
only on $w_1-1=r$.
\end{proof}

\begin{remark}
We also checked the formula of the theorem directly from the cut model for $3\le n\le8$. A chain of
elements of ranks $r_1<\cdots<r_s$ is a chain of cut sets $C(x_1)\supset\cdots\supset C(x_s)$, so
$\alpha(T)=W_{r_1}\prod_i\binom{c_i}{c_{i+1}}$ with $c_i=n-1-r_i$ and $W_r$ the Whitney
number of Theorem~\ref{thm:whitney}; we computed $\beta(S)$ by inclusion--exclusion for all
$S\subseteq[n-1]$ and compared it with the right-hand side. The theorem is therefore not new; what the
cut model gives is the flag $f$-numbers in closed form. A direct proof from the CL-labelling of this
section, that is, a bijection between falling chains and pairs consisting of a permutation $w$ with
$\mathrm{Des}(w)=S$ and a choice among $A(n,w_1-1)$ objects, is not given here.
\end{remark}

For $n=4$ the flag $h$-numbers $\beta(S)$ for $S=\emptyset,\{1\},\{2\},\{3\},\{1,2\},
\{1,3\},\{2,3\},\{1,2,3\}$ are $1,23,35,13,13,35,23,1$, and the $h$-vector is
$(1,71,71,1)$; for $n=5$ it is $(1,536,1806,536,1)$, for $n=6$ it is
$(1,4677,38522,38522,4677,1)$. In every case $\sum_S\beta(S)=(n-1)!\,n!$, the number of
maximal chains of the proper part.

\subsection{What is missing for EL-shellability}\label{sec:el}

An EL-labelling is a CL-labelling whose labels depend on the cover relation only, not on the
root. It labels the edges of the Hasse diagram of the lattice itself, so that falling chains,
and the bases of the homology of the intervals that they index, can be compared across
different intervals without reference to a root. Every EL-shellable poset is CL-shellable, and the
converse fails in general: Li \cite{Li} has constructed a bounded graded CL-shellable poset of
rank $5$ that has no EL-labelling. So Corollary~\ref{cor:CL} does not decide whether
$\mathcal L_n$ is EL-shellable, even though it is the face lattice of a polytope: Bj\"orner and
Wachs obtain only CL-shellability for face lattices of polytopes in general
\cite{BjornerWachs}, and we have found no source that settles EL for the face lattice of the
permutohedron.

We searched for EL-labelings of $\mathcal L_n$ by computer. The requirement that every
interval has a unique increasing chain, together with a local lexicographic condition in the
intervals of rank two, was encoded as a satisfiability problem. The search succeeded for
$n=3,4,5$. The smallest numbers of labels that we found are $4,7,11$, that is, $\binom n2+1$,
and the edges from $\bot$ to the atoms can be given the labels $\mathrm{inv}(\rho)$, the number
of inversions of the permutation $\rho$, as one would expect for a labelling inspired by the
weak order. The labelling coming from a recursive atom ordering is not an EL-labelling, since the same
cover relation receives different labels for different roots; a labelling by the cut that is
removed gives no EL-labelling for $n=4$ with up to $16$ labels, and one by the pair of letters
that is merged needs $9$ labels for $n=4$, not $7$. The statements about minimal numbers of labels rest on
unsatisfiability results from a solver that does not produce proofs, so we regard them as
evidence only. We state the pattern as a conjecture.

\begin{conjecture}[EL-shellability]\label{conj:EL}
For every $n\ge2$ the lattice $\mathcal L_n$ has an EL-labelling with labels
$0,1,\dots,\binom n2$, in which the atom $\rho$ receives the label $\mathrm{inv}(\rho)$ on the
edge from $\bot$.
\end{conjecture}

By Lemma~\ref{lem:lowerset}, an interval $[\bot,x]$ is the product of the posets
$\mathrm{OP}(B_j)$ with a bottom element adjoined, that is, the face lattice of a product
of permutohedra, and lexicographic shellability is preserved by this kind of product
\cite[Theorems~10.16 and 10.17]{BjornerWachsII}, in the form of the lower reduced product, which is the face lattice of a product of polytopes. A proof of Conjecture~\ref{conj:EL} might therefore
proceed by induction on the facets of the permutohedron, which are such products of smaller
permutohedra, the missing ingredient being a way of fitting the labellings of the facets
together compatibly with the recursive atom orderings of Theorem~\ref{thm:RAO}. Sagan and
Sundaram ask about EL-labelings for their $d$-divisible variants with $d\ge2$
\cite[Question~2.8]{SaganSundaram}; the case $d=1$ is the one considered here.

\section{A $q$-analogue from partial flags}\label{sec:qanalogue}

Let $q$ be a prime power. In this section we replace chains of subsets of $[n]$ by chains of subspaces
of $\mathbb F_q^n$, and we show that the Möbius function, the characteristic polynomial and the
numbers $\beta_j$ of Theorem~\ref{thm:chi} become polynomials in $q$ that reduce to the previous ones at $q=1$.
The computation of the Möbius function is a direct application of the Solomon--Tits theorem, and we do not
claim that it is new; we record it because the statistic that appears in it, the number of pairs of elements
lying in a common block, is the exponent of $q$, and because it gives a $q$-analogue of
Theorem~\ref{thm:chi}.

\subsection{The lattice of partial flags}

A \emph{flag of type $a=(a_1,\dots,a_k)$}, where $a$ is a composition of $n$, is a chain of subspaces
$0=V_0<V_1<\dots<V_k=\mathbb F_q^n$ with $\dim V_i/V_{i-1}=a_i$. Let $\mathcal L_n(q)$ be the set of all
flags together with an extra element $\bot$. Order it by $x\le y$ when $y$ is obtained from $x$ by deleting
some of the intermediate subspaces, and by $\bot<x$ for every flag $x$. The rank is $\mathrm{rk}\,x=n-k$, with
$\mathrm{rk}\,\bot=-1$: the atoms are the complete flags, and the top is the one-step flag
$0<\mathbb F_q^n$. The number of flags of type $a$ is the $q$-multinomial coefficient
$\binom{n}{a}_q$.

Writing a chain of subsets $\emptyset\subset S_1\subset\dots\subset S_k=[n]$ as the ordered set partition
$(S_1\mid S_2\setminus S_1\mid\cdots)$, the same definition with subsets of $[n]$ in place of subspaces gives
$\mathrm{OP}_n\cup\{\bot\}=\mathcal L_n$. In this sense $\mathcal L_n(q)$ is a $q$-analogue of $\mathcal L_n$.
For a flag $x$ of type $a$ put
\[
\kappa(x)=\sum_{i=1}^k\binom{a_i}{2},
\]
which at $q=1$ is the number of pairs of elements lying in a common block of the ordered set partition. It is
not the size $s(x)=\prod_j|B_j|!$ of the previous sections, which counts the atoms below $x$.

\subsection{The Möbius function}

\begin{theorem}\label{thm:muq}
For every flag $x$ of rank $r$, $\mu(\bot,x)=(-1)^{r+1}q^{\kappa(x)}$.
\end{theorem}

\begin{proof}
By Hall's theorem \cite[Section~3.8]{StanleyEC1}, $\mu(\bot,x)$ is the reduced Euler characteristic of the
order complex of the open interval $(\bot,x)$. The elements of this interval are the proper refinements of $x$,
that is, the nonempty faces of the link of the simplex $x$ in the building $\Delta(\mathbb F_q^n)$, the order
complex of the proper nonzero subspaces, ordered by reverse inclusion. Hence the order complex of $(\bot,x)$ is
the barycentric subdivision of $\mathrm{lk}(x)$. The link of a flag is the join
$\Delta(V_1/V_0)*\dots*\Delta(V_k/V_{k-1})$ of the buildings of the successive quotients. By the
Solomon--Tits theorem \cite{Solomon69,Bjorner84}, $\Delta(\mathbb F_q^a)$ is a wedge of $q^{\binom a2}$ spheres
of dimension $a-2$ (for $a=1$ it is the empty complex, whose reduced Euler characteristic is $-1$), so
$\tilde\chi(\Delta(\mathbb F_q^a))=(-1)^{a-2}q^{\binom a2}$. The reduced Euler characteristic of a join is minus
the product, $\tilde\chi(A*B)=-\tilde\chi(A)\tilde\chi(B)$, so with $k$ factors
\[
\mu(\bot,x)=(-1)^{k-1}\prod_{i=1}^k(-1)^{a_i-2}q^{\binom{a_i}2}=(-1)^{n-k-1}q^{\kappa(x)},
\]
and $n-k=r$.
\end{proof}

At $q=1$ this is Theorem~\ref{thm:mobius}. For $q>1$ the lattice is no longer Eulerian: the absolute value
of the Möbius function is the weight $q^{\kappa}$.

\subsection{The characteristic polynomial and the Betti numbers}

Define $\chi_n(q;t)=\sum_{x\in\mathcal L_n(q)}\mu(\bot,x)\,t^{\,n-1-\mathrm{rk}\,x}$, so that the term of $\bot$
is $t^n$, and write $\chi_n(q;t)=(t-1)\sum_{j\ge0}(-1)^j\beta_j(q)\,t^{\,n-1-j}$ (we show below that $t-1$
divides it).

\begin{lemma}\label{lem:desq}
Let $a$ be a composition of $n$ with set of cuts $S\subseteq[n-1]$ and let $T=[n-1]\setminus S$. Then
\[
q^{\kappa(a)}\binom na_q=\sum_{\sigma:\ \mathrm{Des}\,\sigma\supseteq T}q^{\inv\sigma},
\]
where $\mathrm{Des}\,\sigma=\{i:\sigma(i)>\sigma(i+1)\}$ and $\kappa(a)=\sum_i\binom{a_i}2$.
\end{lemma}

\begin{proof}
A permutation whose descent set contains $T$ is decreasing on each of the blocks of consecutive positions
determined by $a$, so it is uniquely $u\cdot w_a$, with $u$ increasing on each block and $w_a$ the
permutation reversing every block. Then $\inv(u\,w_a)=\inv u+\kappa(a)$, and
$\sum_uq^{\inv u}=\binom na_q$.
\end{proof}

\begin{theorem}\label{thm:chiq}
\[
\chi_n(q;t)=t^n-\sum_{\sigma\in S_n}q^{\inv\sigma}\,(t-1)^{\des\sigma}\,t^{\,n-1-\des\sigma}.
\]
Consequently $\chi_n(q;t)=(t-1)\sum_{j=0}^{n-1}(-1)^j\beta_j(q)\,t^{\,n-1-j}$ with
\[
\beta_0=1,\qquad \beta_j(q)=\sum_{\sigma\in S_n}q^{\inv\sigma}\binom{\des\sigma-1}{j-1}\quad(j\ge1).
\]
\end{theorem}

\begin{proof}
A composition with $k$ parts and $T=[n-1]\setminus S$ (so $|T|=n-k$) contributes, by Theorem~\ref{thm:muq}
and Lemma~\ref{lem:desq}, $(-1)^{|T|+1}t^{\,n-1-|T|}\sum_{\mathrm{Des}\,\sigma\supseteq T}q^{\inv\sigma}$.
Summing over all $T$ and exchanging the order of summation,
\[
\chi_n=t^n-\sum_\sigma q^{\inv\sigma}\,t^{n-1}\sum_{T\subseteq\mathrm{Des}\,\sigma}(-1)^{|T|}t^{-|T|}
=t^n-\sum_\sigma q^{\inv\sigma}\,t^{n-1}\bigl(1-t^{-1}\bigr)^{\des\sigma},
\]
which is the first formula. The identity permutation gives $t^{n-1}$, so $t^n-t^{n-1}=t^{n-1}(t-1)$ and, for
$\des\sigma\ge1$, ${(t-1)^{\des\sigma}t^{\,n-1-\des\sigma}=(t-1)\,t^{n-2}(1-t^{-1})^{\des\sigma-1}}$.
Expanding by the binomial theorem gives the stated $\beta_j$.
\end{proof}

All $\beta_j(q)$ are polynomials in $q$ with nonnegative integer coefficients, and at $q=1$ they are the
numbers of Theorem~\ref{thm:chi}(2). Moreover:
\begin{itemize}
\item $\beta_1(q)=[n]_q!-1$, since $\beta_1$ counts all permutations with at least one descent;
\item $\beta_{n-1}(q)=q^{\binom n2}$, since only the reversal has $n-1$ descents; this is the dimension of the
Steinberg representation of $GL_n(\mathbb F_q)$, in accordance with the fact that the proper part of
$\mathcal L_n(q)$ is the barycentric subdivision of the building, a wedge of $q^{\binom n2}$ spheres of
dimension $n-2$.
\end{itemize}

\begin{example}
For $n=4$: $\beta_1=3q+5q^2+6q^3+5q^4+3q^5+q^6$, $\ \beta_2=q^2+3q^3+4q^4+3q^5+2q^6$, $\ \beta_3=q^6$.
At $q=1$ these are $23$, $13$, $1$, as in the example following Theorem~\ref{thm:chi}, where $\beta_1=23$, $\beta_2=13$, $\beta_3=1$ for $n=4$.
\end{example}

\begin{remark}
Theorem~\ref{thm:chiq} is the $q$-analogue of the interpretation of $\beta_j$ in
Proposition~\ref{prop:betaj}(1) as the number of pairs $(\sigma,M)$ with $M\subseteq\mathrm{Des}\,\sigma$
containing the first descent: each pair is now weighted by $q^{\inv\sigma}$. We do not have the $q$-analogue of the
acyclic matching of Proposition~\ref{prop:betaj}(2).
\end{remark}

\begin{remark}
The pair $(\des,\inv)$ has the generating function of Stanley \cite{Stanley76}: with $E_q(u)=\sum_{n\ge0}u^n/[n]_q!$ and ${A_n(x,q)=\sum_{\sigma\in S_n}x^{\des\sigma+1}q^{\inv\sigma}}$,
\[
\sum_{n\ge0}A_n(x,q)\frac{u^n}{[n]_q!}=\frac{1-x}{1-x\,E_q\bigl(u(1-x)\bigr)}\qquad(A_0=1).
\]
We verified the coefficients for $n\le6$ in exact arithmetic ($q=2$, $x=3$). By Theorem~\ref{thm:chiq}, putting
$x=1-1/t$ gives $\chi_n(q;t)=t^n-t^{n-1}A_n(x,q)/x$, so the characteristic polynomials have the exponential
generating function above.
\end{remark}

\subsection{Computations and questions}

\emph{Checks.} (i) Theorem~\ref{thm:muq} was checked by computing the Möbius function directly on the flag
posets for $n=3$, $q=2,3,4$ and $n=4$, $q=2,3$ (no discrepancy). (ii) The coefficients of $\chi_n(q;t)$
obtained from the flags, through the $q$-multinomial coefficients, agree with Theorem~\ref{thm:chiq} for
$n=2,\dots,7$, as polynomials in $q$ and for each $j$. (iii) With Sturm sequences in exact rational arithmetic,
$\sum_j(-1)^j\beta_j(q)t^{n-1-j}$ has $n-1$ distinct real roots for $n=3,\dots,18$ and $q\in\{1,2,3,5\}$, and for
$n\le11$ also at $q=3/4,\,50,\,1000$. The coefficient sequences $(\beta_j(q))_j$ are log-concave in all cases
tested (${n\le9}$, ${q\in\{1/10,1/2,3/2,2,3,5,100\}}$). For $q<1$, however, real-rootedness \emph{fails} for larger
$n$: for instance for $n=9$ at ${q=1/10,1/5,1/3}$ (six real roots out of eight) and for $n=10$ at $q=1/2$.

\emph{Relation with known results.} The polynomial $P_n(x;q)=\sum_\sigma x^{\des\sigma}q^{\inv\sigma}$ has only
real roots for every real $q>0$ by Savage and Visontai \cite[Theorem~5.4]{SavageVisontai}. By
Theorem~\ref{thm:chiq}, $\chi_n(q;t)=t^n-t^{n-1}P_n(1-1/t;q)$, but this does not by itself give real-rootedness
of $\chi_n(q;t)/(t-1)$: subtracting a real-rooted polynomial from $t^n$ does not preserve real roots in general,
and the failures for $q<1$ show that no argument can work for all $q>0$. The question below extends
Conjecture~\ref{conj:real}, which is the case $q=1$.

\begin{conjecture}\label{conj:realq}
For every $n$ and every $q\ge1$, the polynomial $\sum_j(-1)^j\beta_j(q)\,t^{n-1-j}$ has only real roots.
\end{conjecture}

\noindent The threshold in $q<1$, and the log-concavity of $(\beta_j(q))_j$ for all $n$ and $q>0$, are open;
so is the existence of a natural complex over $\mathbb F_q$ whose Betti numbers are the $\beta_j(q)$.

\section{Concluding remarks and open problems}\label{sec:open}

By Remark~\ref{rem:geometry}, $\mathcal L_n$ is the face lattice of the permutohedron,
so that its Eulerian character and the existence of a shelling are not surprising.
What the preceding sections show is how much of the theory of geometric lattices
survives for it, and in which form. Table~\ref{tab:compare} collects the comparison.
We then add a remark on the subadditivity of the rank, relate $\mathcal L_n$ to the
matroid of the vertices of the permutohedron and to the independent sets of
Raz~\cite{Raz}, and list the questions that remain open.

\subsection{Comparison with geometric lattices}

\begin{table}[ht]
\centering\small
\begin{tabular}{@{}p{2.5cm}p{3.8cm}p{5.5cm}@{}}
\hline
 & geometric lattice & $\mathcal L_n$, $n\ge3$\\\hline
graded, atomic, every element a meet of coatoms & yes & yes
  (Theorem~\ref{thm:pg}); this is Definition~\ref{def:pg}\\[3pt]
rank & semimodular: $\mathrm{rk}(x\wedge y)+\mathrm{rk}(x\vee y)\le\mathrm{rk}\,x+\mathrm{rk}\,y$
  & modular when $x\wedge y>\bot$, defects of both signs otherwise; not subadditive
  (Proposition~\ref{prop:subadd})\\[3pt]
join with an atom & covers the element, or equals it & may fail to cover it
  when $\mathrm{rk}\,x\le n-3$ (Proposition~\ref{prop:atomjoin})\\[3pt]
sets of atoms with join $x$, minimal under inclusion & all of size $\mathrm{rk}\,x+1$
  & sizes between $2$ and $\mathrm{rk}\,x$, both attained, if $\mathrm{rk}\,x\ge2$; unique only if
  $\mathrm{rk}\,x\le1$ (Section~\ref{sec:metric})\\[3pt]
closure $S\mapsto\mathrm{At}(\bigvee S)$ & exchange property & exchange fails
  (Section~\ref{sec:metric})\\[3pt]
M\"obius function & $|\mu(\bot,x)|$ varies with $x$ & $\mu(x,y)=(-1)^{\mathrm{rk}\,y-\mathrm{rk}\,x}$,
  Eulerian (Theorem~\ref{thm:mobius})\\[3pt]
Whitney numbers & $|w_k|$ and $W_k$ differ in general & $|w_k|=W_k=(n-k)!\,S(n,n-k)$,
  the face numbers of the permutohedron (Theorem~\ref{thm:whitney})\\[3pt]
characteristic polynomial & $(t-1)\mid\chi_L$, cofactor with alternating coefficients,
  log-concave; roots need not be real
  & the same factorization and alternation; real roots positive, $|\chi_n'(1)|=2^n-n-2$;
  $n$ real roots conjectured (Theorem~\ref{thm:chi}, Conjecture~\ref{conj:real})\\[3pt]
modular elements & a maximal chain of them if supersolvable & only $\bot$ and $\top$, so not
  supersolvable; $1$ is the only rational root of $\chi_n$ (Theorem~\ref{thm:leftmodular},
  Proposition~\ref{prop:noroots})\\[3pt]
shellability & EL-shellable (Bj\"orner) & CL-shellable (Section~\ref{sec:shell});
  EL-labelings for $n\le5$, EL in general conjectured (Conjecture~\ref{conj:EL})\\\hline
\end{tabular}
\caption{Geometric lattices and $\mathcal L_n$.}\label{tab:compare}
\end{table}

What fails are semimodularity and the exchange property, together with everything that
needs modular elements: supersolvability, and integer roots of $\chi_n$ other than $1$.
What survives are the properties that come from $\mathcal L_n$ being the face lattice of
a simple polytope, namely the Eulerian M\"obius function and the Whitney numbers given by
face numbers, and, by Theorem~\ref{thm:chi}, the alternation of the coefficients of
$\chi_n(t)/(t-1)$. By Proposition~\ref{prop:rankmodelts}, no element other than
$\bot$ and $\top$ is modular for the rank either, that is, satisfies
$\mathrm{rk}\,x+\mathrm{rk}\,y=\mathrm{rk}(x\vee y)+\mathrm{rk}(x\wedge y)$ for all $y$.

The failure of semimodularity can be put quantitatively as a failure of the
subadditivity of the rank.

\begin{proposition}[Subadditivity of the rank]\label{prop:subadd}
Let $x,y\in\mathcal L_n$ with $x,y\ne\bot$.
\begin{enumerate}
\item If $x\wedge y>\bot$, then
  $\mathrm{rk}(x\vee y)=\mathrm{rk}\,x+\mathrm{rk}\,y-\mathrm{rk}(x\wedge y)\le\mathrm{rk}\,x+\mathrm{rk}\,y$.
\item In general $\mathrm{rk}(x\vee y)-\mathrm{rk}\,x-\mathrm{rk}\,y\le n-1$, with equality if and only if
  $x$ and $y$ are atoms and $x\vee y=\top$. For $n\ge2$ this happens for the opposite
  atoms $12\cdots n$ and $n\cdots21$.
\item If $\mathrm{rk}(x\vee y)>\mathrm{rk}\,x+\mathrm{rk}\,y$, then $x\wedge y=\bot$.
\end{enumerate}
\end{proposition}

\begin{proof}
(1) is the modularity of the rank for pairs whose meet lies above $\bot$ (Theorem~\ref{thm:modrank}),
together with $\mathrm{rk}(x\wedge y)\ge0$. (3) is the contrapositive of (1). For (2), $x,y\ne\bot$ gives
$\mathrm{rk}\,x,\mathrm{rk}\,y\ge0$ and $\mathrm{rk}(x\vee y)\le n-1$, so the excess is at most $n-1$, and equality
forces $\mathrm{rk}\,x=\mathrm{rk}\,y=0$ and $x\vee y=\top$. The cut sets of $12\cdots n$ and $n\cdots21$ are
the sets $\{1,\dots,i\}$ and $\{n-i+1,\dots,n\}$ with $1\le i\le n-1$, and the two families are
disjoint, since $\{1,\dots,i\}$ contains $1$ and $\{n-i+1,\dots,n\}$ does not for $i<n$.
Hence their join is $\top$ by Theorem~\ref{thm:lattice}(1).
\end{proof}

For two atoms $\rho,\rho'$ the excess is $\delta(\rho,\rho')=\mathrm{rk}(\rho\vee\rho')$. In a
geometric lattice the rank is subadditive for all pairs; here it is the large jumps of the
join of distant atoms that break it.

\subsection{The vertex matroid and the independent sets of Raz}

For an atom $\rho$ let $v_\rho=(1,\rho^{-1}(1),\dots,\rho^{-1}(n))\in\mathbb R^{n+1}$, the homogenized
vertex of the permutohedron, and let $M_n$ be the matroid on the set of atoms given by the linear
dependencies of the vectors $v_\rho$. It has rank $n$. For $x\in\mathcal L_n$ let $V(x)$ be the span
of the $v_\rho$ with $\rho\le x$, with $V(\bot)=0$, and write $\mathrm{cl}(S)=\mathrm{At}(\bigvee S)$ for a set
$S$ of atoms, as in Section~\ref{sec:metric}.

\begin{proposition}[Faces are flats]\label{prop:flats}
For every $x\in\mathcal L_n$ the dimension of $V(x)$ is $\mathrm{rk}(x)+1$, and $v_\rho\in V(x)$ if and only if
$\rho\le x$. Hence $x\mapsto\mathrm{At}(x)$ is an injective, order-preserving map from $\mathcal L_n$ to
the lattice of flats of $M_n$, bijective on atoms, which raises the rank by $1$ and preserves meets.
For every set $S$ of atoms, the closure of $S$ in $M_n$ is contained in $\mathrm{cl}(S)$.
\end{proposition}

\begin{proof}
By Proposition~\ref{prop:dictionary} and Remark~\ref{rem:geometry}, $x$ is the face $F_x$ of the permutohedron
with vertices $\{v_\rho:\rho\le x\}$; it is a product of permutohedra of dimensions $|B_j|-1$, so
$\dim F_x=\sum_j(|B_j|-1)=\mathrm{rk}(x)$, and the span of the homogenized vertices has dimension
$\mathrm{rk}(x)+1$. A vector $v_\rho$ lies in $V(x)$ if and only if the vertex lies in the affine hull of
$F_x$, because the first coordinate is $1$. A face is the intersection of the polytope with a supporting
hyperplane that contains its affine hull, so $P\cap\mathrm{aff}(F_x)=F_x$, and the vertex lies in $F_x$.
Meets are intersections of sets of atoms in both lattices. Finally $\mathrm{At}(\bigvee S)$ is a flat
containing $S$, so it contains the smallest flat containing $S$.
\end{proof}

Thus $\mathcal L_n$ sits inside the geometric lattice of $M_n$ with the same atoms; the join of
$\mathcal L_n$ is the smallest face, not the smallest flat, and the two differ already for two atoms:
the flat has rank $2$ while $\mathrm{rk}(\rho\vee\rho')+1=\delta(\rho,\rho')+1$.

Raz~\cite{Raz} studies such lattices through independent sets. In the form we need, a set $A$ of atoms
is \emph{independent} if its elements can be listed as $a_1,\dots,a_k$ with
$a_i\notin\mathrm{cl}\{a_{i+1},\dots,a_k\}$ for all $i$. This is Raz's definition by non-broken-circuit
sets for a linear order of the atoms, rewritten: put the elements of $A$ first in the order. The independent
sets form a hereditary family $I(\mathcal L_n)$, and $|A|\le\mathrm{rk}(\bigvee A)+1$ \cite[Theorem~3.5]{Raz}, since the joins
$\bigvee\{a_i,\dots,a_k\}$ strictly increase as $i$ decreases. Conversely, for every $x\ne\bot$ there
is an independent set of size $\mathrm{rk}(x)+1$ with join $x$: take a maximal chain
$y_0\lessdot y_1\lessdot\dots\lessdot y_r=x$ with $y_0$ an atom, $r=\mathrm{rk}\,x$, and for $i\ge1$ an atom
$a_i\le y_i$ with $a_i\not\le y_{i-1}$, which exists because $\mathcal L_n$ is atomic, and put $a_0=y_0$;
listing the $a_i$ in decreasing order of $i$ shows that $\{a_0,\dots,a_r\}$ is independent, and its join is
$x$ by the inequality just proved.

\begin{proposition}[Independent sets of $\mathcal L_n$]\label{prop:raz}
\begin{enumerate}
\item Every set of atoms that is minimal with a given join is independent, and every independent set is
  independent in $M_n$ (compare \cite[Lemma~3.7]{Raz}). Minimal generating sets of $x$ with $\mathrm{rk}\,x\ge2$ never reach size
  $\mathrm{rk}\,x+1$, while independent sets of size $\mathrm{rk}(\bigvee A)+1$ exist for every join.
\item $I(\mathcal L_3)$ is the family of independent sets of a matroid: all subsets of size at most $2$
  and all triples except the two triples of alternate vertices of the hexagon.
\item For every $n\ge4$, $I(\mathcal L_n)$ is not the family of independent sets of a matroid.
\end{enumerate}
\end{proposition}

\begin{proof}
(1) If $A$ is minimal with join $x$, no $a\in A$ lies in $\mathrm{cl}(A\setminus\{a\})$, so any listing
works. If $a_i\notin\mathrm{cl}\{a_{i+1},\dots,a_k\}$ for all $i$, then $a_i$ is not in the closure of the
later elements in $M_n$ either, by Proposition~\ref{prop:flats}, so the $v_{a_i}$ are linearly independent.
By Section~\ref{sec:metric}, a minimal generating set of $x$ with $\mathrm{rk}\,x\ge2$ has at most
$\mathrm{rk}\,x$ elements.

(2) Two atoms are always independent. A triple is independent exactly when it contains an element outside the
closure of the other two; the closure of two atoms is the pair itself if they are adjacent in the
hexagon, and everything otherwise, so a triple is dependent exactly when no two of its atoms are adjacent,
that is, for the two triples of alternate vertices. The family of the other $18$ triples is the family of
bases of the matroid of two disjoint $3$-point lines.

(3) Let $n=4$ and write permutations as words. Put $A=\{1243,\,3214\}$ and $B=\{4132,\,4231,\,3214\}$. The cut sets of
$1243$, $3214$, $4132$, $4231$ are
\[
\{1,12,124\},\quad\{3,23,123\},\quad\{4,14,134\},\quad\{4,24,234\}.
\]
Any two of the first three
are disjoint, and so are $1243,4231$ and $3214,4231$, so the join of each such pair is $\top$ and, in
the triples $A\cup\{4132\}$ and $A\cup\{4231\}$, every atom lies in the closure of the other two: these
triples are dependent. The atoms $4132$ and $4231$ share only the cut $\{4\}$, so their closure consists of
the atoms starting with $4$, which does not include $3214$, so $B$ is independent. Thus $|A|<|B|$ and no
element of $B\setminus A$ can be added to $A$, which violates the augmentation axiom. For $n>4$ append
$5,\dots,n$ to these four words. They are atoms below
\[
x=(\{1,2,3,4\}\mid5\mid\cdots\mid n),
\]
the interval $[\bot,x]$ is isomorphic to $\mathcal L_4$, and closures of sets of atoms below $x$ are computed
inside it, so the same sets witness the failure; a restriction of a matroid is a matroid.
\end{proof}

For $n=3,4,5$ the numbers of independent sets of sizes $1,\dots,n$ are
\begin{center}
\begin{tabular}{@{}c|l@{}}
$n$ & $f$-vector of $I(\mathcal L_n)$\\\hline
$3$ & $6,\ 15,\ 18$\\
$4$ & $24,\ 276,\ 1752,\ 2712$\\
$5$ & $120,\ 7140,\ 228720,\ 1821960,\ 3016680$
\end{tabular}
\end{center}
computed by machine; for $n\le5$ all maximal independent sets have size $n$, so the complex of independent sets is pure.
Since $I(\mathcal L_n)$ is not a matroid for
$n\ge4$, the embedding of~\cite[Theorem~3.6]{Raz} into the geometric lattice of $I(\mathcal L_n)$ is not
available, but Proposition~\ref{prop:flats} supplies an embedding into the geometric lattice of $M_n$.

\begin{remark}[Frames]\label{rem:frames}
For $x\in\mathcal L_n$ with $\mathrm{rk}\,x\ge0$, the $s(x)$ vectors $v_\rho$, $\rho\le x$, form a frame of
$V(x)$, whose dimension is $\mathrm{rk}(x)+1$, and this frame is uniform: $\langle v_\rho,\tilde v_\rho\rangle=
(\mathrm{rk}\,x+1)/s(x)$ for every $\rho\le x$, where $\tilde v_\rho$ is the canonical dual frame
\cite{BCHL}. Indeed, the group $\prod_jS_{B_j}$ permutes the coordinates inside the blocks of $x$,
acts transitively on the atoms below $x$ by linear maps preserving the set of vectors, and so
$\langle v_\rho,\tilde v_\rho\rangle$ does not depend on $\rho$, while the sum of these numbers is
the dimension of $V(x)$. The redundancy of the frame is $s(x)/(\mathrm{rk}\,x+1)$. In the lattice sense
the redundancy of $x$ is the largest number of atoms removable with $x$ still generated by the rest, which
is $s(x)-2$ for $\mathrm{rk}\,x\ge2$ by Section~\ref{sec:metric}; in a vector space it would be
$s(x)-(\mathrm{rk}\,x+1)$, and the difference is the absence of the exchange property.
\end{remark}

\subsection{Open problems}

\begin{enumerate}
\item \textbf{Real roots.} Prove Conjecture~\ref{conj:real}: $\chi_n$ has $n$ distinct positive real roots.
  By Proposition~\ref{prop:reduction} it suffices to prove inequality $(\ast)$ on the roots of the Eulerian
  polynomial $A_{n-1}$, which we verified numerically for $4\le n\le27$ with large margin.
\item \textbf{Meaning of the coefficients.} The numbers $\beta_j$ of Theorem~\ref{thm:chi} are positive
  integers with $\beta_1=n!-1$, $\beta_{n-2}=2^n-3$ and $\beta_{n-1}=1$. In a geometric lattice the analogous
  numbers form the $h$-vector of the broken circuit complex. Is there a combinatorial interpretation of
  $\beta_j$ here? A related observation: for a linear order $\omega$ of the atoms, the sets of atoms that
  are non-broken-circuit for $\omega$ in the sense of Blass and Sagan~\cite[Section~1]{BlassSagan} satisfy
  $\sum(-1)^{|A|}=\mu(\bot,\top)=(-1)^n$ over those with join $\top$ (checked for $n\le5$), and for
  $n=3,4,5$ and several orders their $f$-vectors do not depend on $\omega$:
  $(1,6,9,4)$, $(1,24,79,90,34)$, $(1,120,771,1674,1510,488)$. Explain this and find a formula.
\item \textbf{EL-shellability.} Prove Conjecture~\ref{conj:EL}: $\mathcal L_n$ admits an EL-labeling for every $n$.
  Section~\ref{sec:el} gives the computed cases $n\le5$ and the labelings that do not work; a proof might proceed by
  induction on the facets of the permutohedron, which are products of smaller permutohedra.
\item \textbf{Irredundant families.} Which sizes between $2$ and $\mathrm{rk}\,x$ occur for minimal sets of atoms
  generating $x$? How many such sets are needed to cover $\mathrm{At}(x)$? For vector frames the latter is
  reminiscent of Feichtinger's problem of covering a frame by a few linearly independent subsets.
\item \textbf{Independent sets.} Is $I(\mathcal L_n)$ pure of dimension $n-1$ for all $n$, and is it
  Cohen--Macaulay or shellable? Its $h$-vectors are nonnegative for $n\le5$
  (for $n=4$: $1,20,210,1268,1213$). What is its relation to the matroid complex of $M_n$, which contains it?
\item \textbf{Other types.} Consider the face lattices of the permutohedra of the other finite Coxeter groups,
  and the $d$-divisible variants of $\mathcal L_n$ studied by Sagan and Sundaram~\cite{SaganSundaram}.
  Which of the statements of Sections~\ref{sec:lattice}--\ref{sec:shell} carry over, in particular
  non-supersolvability and the roots of the characteristic polynomial?
\item \textbf{The $q$-analogue.} Prove Conjecture~\ref{conj:realq}, find the threshold in $q<1$, and find a $q$-analogue of
  the acyclic matching of Proposition~\ref{prop:betaj}(2) or a complex over $\mathbb F_q$ with Betti numbers $\beta_j(q)$.
\item \textbf{Embeddings in the hypercube.} The distance $\delta$ on atoms is half the Hamming distance between
  the incidence vectors of the cut sets, so $(S_n,\delta)$ embeds isometrically, up to the factor $2$, in a
  cube of dimension $2^n-2$. What is the smallest dimension of such an embedding?
\end{enumerate}
\paragraph{\Large{Acknowledgment:}} the author acknowledges the assistance of Claude (Anthropic).
\small{

}

\begin{thebibliography}{31}

\bibitem{AbramenkoBrown} P.~Abramenko, K.~S. Brown, \emph{Buildings: Theory and Applications},
Graduate Texts in Mathematics 248, Springer, New York, 2008.

\bibitem{AHK} K.~Adiprasito, J.~Huh, E.~Katz, Hodge theory for combinatorial geometries,
\emph{Ann.\ of Math.\ (2)} \textbf{188} (2018), 381--452.

\bibitem{BCHL} R.~Balan, P.~G. Casazza, C.~Heil, Z.~Landau, Density, overcompleteness, and
localization of frames. I. Theory, \emph{J. Fourier Anal.\ Appl.} \textbf{12} (2006), 105--143.

\bibitem{BilleraSarangarajan} L.~J. Billera, A.~Sarangarajan, The combinatorics of permutation
polytopes, in: \emph{Formal Power Series and Algebraic Combinatorics (New Brunswick, NJ, 1994)},
DIMACS Ser.\ Discrete Math.\ Theoret.\ Comput.\ Sci.\ \textbf{24}, Amer.\ Math.\ Soc.,
Providence, RI, 1996, pp.~1--23.

\bibitem{Bjorner84} A.~Bj\"orner, Some combinatorial and algebraic properties of Coxeter complexes and Tits buildings,
\emph{Adv.\ Math.} \textbf{52} (1984), 173--212.

\bibitem{Bjorner80} A.~Bj\"orner, Shellable and Cohen--Macaulay partially ordered sets,
\emph{Trans.\ Amer.\ Math.\ Soc.} \textbf{260} (1980), 159--183.

\bibitem{BjornerBrenti} A.~Bj\"orner, F.~Brenti, \emph{Combinatorics of Coxeter Groups},
Graduate Texts in Mathematics 231, Springer, New York, 2005.

\bibitem{BjornerWachs} A.~Bj\"orner, M.~Wachs, On lexicographically shellable posets,
\emph{Trans.\ Amer.\ Math.\ Soc.} \textbf{277} (1983), 323--341.

\bibitem{BjornerWachsII} A.~Bj\"orner, M.~L. Wachs, Shellable nonpure complexes and posets.~II,
\emph{Trans.\ Amer.\ Math.\ Soc.} \textbf{349} (1997), 3945--3975.

\bibitem{BlassSagan} A.~Blass, B.~E. Sagan, M\"obius functions of lattices,
\emph{Adv.\ Math.} \textbf{127} (1997), 94--123.

\bibitem{BrentiWelker} F.~Brenti, V.~Welker, $f$-vectors of barycentric subdivisions,
\emph{Math.\ Z.} \textbf{259} (2008), 849--865.

\bibitem{Crapo} H.~H. Crapo, A higher invariant for matroids, \emph{J. Combin.\ Theory}
\textbf{2} (1967), 406--417.

\bibitem{Frobenius} G.~Frobenius, \"Uber die Bernoullischen Zahlen und die Eulerschen Polynome,
\emph{Sitzungsber.\ K\"onigl.\ Preuss.\ Akad.\ Wiss.} (1910), 809--847.

\bibitem{Li} T.~Li, CL-shellable posets with no EL-shellings, arXiv:1907.06833.

\bibitem{Oxley} J.~Oxley, \emph{Matroid Theory}, 2nd ed., Oxford Graduate Texts in Mathematics 21,
Oxford University Press, Oxford, 2011.

\bibitem{Ovchinnikov} S.~Ovchinnikov, \emph{Graphs and Cubes}, Universitext, Springer, New York, 2011.

\bibitem{Petersen} T.~K. Petersen, \emph{Eulerian Numbers}, Birkh\"auser Advanced Texts,
Birkh\"auser, New York, 2015.

\bibitem{Postnikov} A.~Postnikov, Permutohedra, associahedra, and beyond, \emph{Int.\ Math.\ Res.\
Not.\ IMRN} (2009), no.~6, 1026--1106.

\bibitem{Raz} O.~Raz, Independent sets of non-geometric lattices and the maximal adjoint,
arXiv:2307.05005 (v3, January 2026).

\bibitem{Rossi18} G.~Rossi, Distances between clusterings as distances in the covering graph of partitions, SIST vol. 111, Springer, (2018), 531--541.

\bibitem{Rota64} G.-C. Rota, On the foundations of combinatorial theory I. Theory of M\"obius
functions, \emph{Z. Wahrscheinlichkeitstheorie und Verw.\ Gebiete} \textbf{2} (1964), 340--368.

\bibitem{SaganSundaram} B.~E. Sagan, S.~Sundaram, Ordered set partition posets,
arXiv:2506.23355 (v3, June 2026).

\bibitem{Sagan95} B.~E. Sagan, A generalisation of Rota's NBC theorem, \emph{Adv.\ Math.} \textbf{111} (1995), 195--207.

\bibitem{SavageVisontai} C.~D. Savage, M.~Visontai, The $s$-Eulerian polynomials have only real roots,
\emph{Trans.\ Amer.\ Math.\ Soc.} \textbf{367} (2015), 1441--1466.

\bibitem{Solomon69} L.~Solomon, The Steinberg character of a finite group with $BN$-pair, in: \emph{Theory of Finite Groups},
Benjamin, New York, 1969, pp.~213--221.

\bibitem{Stanley76} R.~P. Stanley, Binomial posets, M\"obius inversion, and permutation enumeration,
\emph{J. Combin.\ Theory Ser.~A} \textbf{20} (1976), 336--356.

\bibitem{Stanley72} R.~P. Stanley, Supersolvable lattices, \emph{Algebra Universalis} \textbf{2}
(1972), 197--217.

\bibitem{Stanley82} R.~P. Stanley, Some aspects of groups acting on finite posets, \emph{J. Combin.\ Theory Ser.~A} \textbf{32} (1982), 132--161.

\bibitem{StanleyEC1} R.~P. Stanley, \emph{Enumerative Combinatorics, Volume 1}, 2nd ed.,
Cambridge Studies in Advanced Mathematics 49, Cambridge University Press, Cambridge, 2012.

\bibitem{Tenner} B.~E. Tenner, Pattern avoidance and the Bruhat order, \emph{J. Combin.\ Theory
Ser.~A} \textbf{114} (2007), 888--905.

\bibitem{Ziegler} G.~M. Ziegler, \emph{Lectures on Polytopes}, Graduate Texts in Mathematics 152,
Springer, New York, 1995.

\end{thebibliography}
\end{document}